\documentclass[12pt]{article}

\usepackage{amssymb,amsmath,amsfonts,bm}
\usepackage[margin=2cm]{geometry}
\usepackage{cuted}
\usepackage[ansinew]{inputenc}
\usepackage[english]{babel}
\usepackage{fancyhdr}
\usepackage{anyfontsize}
\usepackage[active]{srcltx}
\usepackage{soul}
\usepackage{xcolor}
\usepackage[colorlinks=true]{hyperref}
\usepackage{color}
\usepackage{array}
\usepackage{textcomp}
\usepackage[linesnumbered,ruled,vlined]{algorithm2e}
\usepackage{pgfplots}
\pgfplotsset{compat=1.17}
\newcommand{\matrices}[1]{\mathbf{#1}}
\newcommand{\functions}[1]{\bm{#1}}
\newcommand{\vectors}[1]{\mathbf{#1}}

\renewcommand{\P}[2][p]{{\left[#2\right]_{#1}}}
\newcommand{\alphavec}{\boldsymbol \alpha}

\newcommand{\mus}{\mu\text{s}}
\newcommand{\mum}{\mu\text{m}}

\newcommand{\mR}{\mathcal{R}}

\newcommand{\y}{\vectors{y}}

\newcommand{\z}{\vectors{z}}
\newcommand{\C}{\vectors{C}}
\newcommand{\0}{\vectors{0}}
\newcommand{\W}{\vectors{W}}

\newcommand{\f}{\vectors{f}}

\newcommand{\X}{\vectors{X}}
\newcommand{\Y}{\vectors{Y}}

\newcommand{\R}{\vectors{R}}
\newcommand{\Svec}{\vectors{S}}

\newcommand{\N}{\vectors{N}}

\newcommand{\B}{\matrices{B}}
\newcommand{\A}{\matrices{A}}
\newcommand{\I}{\matrices{I}}
\newcommand{\Q}{\functions{Q}}

\newcommand{\M}{\matrices{M}}
\newcommand{\K}{\matrices{K}}
\newcommand{\U}{\matrices{U}}
\newcommand{\V}{\matrices{V}}
\newcommand{\G}{\matrices{G}}
\newcommand{\Hv}{\matrices{H}}
\newcommand{\F}{\matrices{F}}
\newcommand{\D}{\matrices{D}}
\newcommand{\E}{\matrices{E}}
\newcommand{\phivec}{\boldsymbol \phi}

\newcommand{\Lambdavec}{\boldsymbol \Lambda}
\newcommand{\Xivec}{\boldsymbol \Xi}

\newcommand{\xivec}{\boldsymbol \xi}
\newcommand{\muvec}{\boldsymbol \mu}
\newcommand{\nuvec}{\boldsymbol \nu}

\newcommand{\Psivec}{\boldsymbol \Psi}
\newcommand{\Upsvec}{\boldsymbol \Upsilon}

\newcommand{\dblnk}{\begin{minipage}{12pt}\fontsize{6pt}{3pt}\selectfont{$>1$\\$<p$}\end{minipage}}

\DeclareMathOperator{\e}{e}
\usepackage{subfiles}
\usepackage{subfig}
\usepackage{smartdiagram}
\usesmartdiagramlibrary{additions}
\usepackage{lineno}

\usepackage{upgreek}
\usepackage{textgreek}

\usepackage{cancel}

\usepackage{authblk}

\usepackage[symbol]{footmisc}

\providecommand{\keywords}[1]{\textbf{\textit{keywords---}} #1}

\title{Direct parametrisation of invariant manifolds for non-autonomous forced systems including superharmonic resonances}

\author{Alessandra Vizzaccaro$^1$, Giorgio Gobat$^2$\footnote{corresponding author - giorgio.gobat@polimi.it}, Attilio Frangi$^2$, Cyril Touz{\'e}$^3$}

\date{{\small
    $^1$ College of Engineering, Mathematics and Physical Sciences, University of Exeter, Exeter, UK\\%
    $^2$ Department of Civil and Environmental Engineering, Politecnico di Milano, Milano, Italy\\%
    $^3$ Institute of Mechanical Sciences and Industrial Applications (IMSIA), ENSTA Paris - CNRS - EDF - CEA, Institut Polytechnique de Paris, Palaiseau, France\\[2ex]}
}
\begin{document}
\maketitle

\begin{abstract}
The direct parametrisation method for invariant manifold is a model-order reduction technique that can be applied to nonlinear systems described by PDEs and discretised {\it e.g.} with a finite element procedure in order to derive efficient reduced-order models (ROMs).
In nonlinear vibrations, it has already been applied to autonomous and non-autonomous problems to propose ROMs that can compute backbone and frequency-response curves of structures with geometric nonlinearity. While previous developments used a first-order expansion to cope with the non-autonomous term, this assumption is here relaxed by proposing a different treatment. The key idea is to enlarge the dimension of the parametrising coordinates with additional entries related to the forcing.
A new algorithm is derived with this starting assumption and, as a key consequence, the resonance relationships appearing through the homological equations involve multiple occurrences of the forcing frequency, showing that with this new development, ROMs for systems exhibiting a superharmonic resonance, can be derived. The method is implemented and validated on academic test cases involving beams and arches. It is numerically demonstrated that the method generates efficient ROMs for problems involving 3:1 and 2:1 superharmonic resonances, as well as converged results for systems where the first-order truncation on the non-autonomous term showed a clear limitation.
\end{abstract}


\section{Introduction}

Model order reduction for nonlinear vibratory systems using the nonlinear normal modes (NNMs) of the systems, has been under development since the first works by Rosenberg in the 1960s~\cite{Rosenberg62}. In the 1990s, the concept regained attention and the link with Lyapunov subcentre manifold, and more generally speaking invariant manifolds of the system, was first introduced~\cite{VakakisNNM,ShawPierre91,ShawPierre93,touze03-NNM}. 
Despite efficient results obtained in early 2000, see {\it e.g.}~\cite{PesheckJSV,JIANG2005H,TOUZE:CMAME:2008,LegrandIM}, the method has not been broadly adopted, probably because direct applications to finite element (FE) problems have remained elusive. 
However, recent developments overcame this issue by proposing two important advances allowing one to set the theoretical settings in a unified framework, together with proposing direct computations, applicable to structures discretised by the FE approach.

The first important step ahead has been to unify the different computational approaches thanks to the parametrisation method for invariant manifold, a technique first proposed in the dynamical systems community~\cite{Cabre1,Cabre3,Haro}. Rewriting the problem from the invariance equation, one is then able to link the different computations using either centre manifold technique or normal form approach as different parametrisation styles~\cite{Haro,ReviewROMGEOMNL}. 
In the dynamical systems community, the parametrisation method has been first applied to different limit sets to extend the analysis to {\it e.g.}  limit cycles or stable manifolds~\cite{Castelli2015,vdbJay2016}. 
Applications to problems semi-discretised with FE have been proposed~\cite{Gonzalez2022}, whereas applications to celestial mechanics have been more largely developed, see {\it e.g.}~\cite{Haro,LeBihan_2017} and references therein. 
Recent developments involve for example interesting discussions on the validity range of the different styles of parametrisations~\cite{Stoychev:failing}, and application of the parametrisation method to the Navier-Stokes equation~\cite{Buza:NS}.
In the field of nonlinear vibrations, it has been first used in~\cite{Haller2016,PONSIOEN2018} for model order reduction.
However, the calculations used, as a first step, the equations of motion expressed in the modal basis, which was a known limitation for applications of invariant manifold methods to large-scale FE problems~\cite{TOUZE:JSV:2006,ReviewROMGEOMNL}. 
In conjunction with this technical development, a key result shown in~\cite{Haller2016} relates to the existence and uniqueness of invariant manifolds used for model-order reduction, which has been named as Spectral Submanifold (SSM) for damped systems. Thanks to this proof, an effective mathematical framework has been settled, which gives a ground existence to the nonlinear normal modes (NNMs) defined as invariant manifolds. This seminal result allows understanding that the SSM is the targeted unique invariant structure in phase space that can be used as a reduced-order model. Although in a conservative framework an NNM is largely understood as being the associated Lyapunov subcenter manifold (LSM), in the dissipative case an NNM must be defined as the associated SSM. 
Retrospectively, both earlier works, as those reported in~\cite{ShawPierre93,PesheckJSV,TOUZE:CMAME:2008} and more recent ones, shown for instance in~\cite{Haller2016,artDNF2020,vizza21high,JAIN2021How}, compute numerical approximations of SSMs and LSMs.

The second important step has been to propose a direct approach such that one could bypass the projection onto the linear modes basis as a starting point. 
Indeed, in a FE context, this step is a clear and known limitation that completely prevents a broad application of invariant-based methods to engineering structures. 
Direct calculations have been proposed in~\cite{artDNF2020}, where the normal form approach derived in~\cite{touze03-NNM,TOUZE:JSV:2006} has been adapted such that the calculations can be operated in a direct and non-intrusive manner. 
The method has then been used in a more general setting, including internal resonance and engineering applications~\cite{AndreaROM}.
Then, direct computations relying on the parametrisation method, have been proposed in~\cite{JAIN2021How,li2021periodic} and, by the authors, in~\cite{vizza21high,opreni22high}. In conjunction with these developments, efficient open-source softwares implementing the developments have been released in order to share the method: \texttt{SSMtools}~\cite{ssmtools} and \texttt{MORFE}~\cite{morfe}. Note that the present authors adopted the acronym DPIM (direct parametrisation of invariant manifold) in their works to refer to the algorithmic procedure required to compute approximations of various objects: SSMs in the dissipative case, LSMs in the conservative case, and whiskers in the non-autonomous case~\cite{haro2006parameterization2,opreni22high,JAIN2021How}.

These implementations of the parametrisation method have then been applied to various problems~\cite{li2021periodic,li2021bifurcation}, including rotating structures~\cite{Martin:rotation} and piezo-electrically actuated MEMS~\cite{opreniPiezo}.
Comparisons to deep-learning-based approaches are also reported in~\cite{Gobat:DLDPIM}, underlining how well-tuned data-driven techniques can  recover the invariant manifold solution.

Taking into account an external forcing makes the dynamical system non-autonomous, such that the invariant manifolds become time-dependent, as opposed to those of the unforced system. In dynamical systems literature, early developments considering the forcing are available for centre manifold approach~\cite{cox_roberts_1991}, or normal form theory~\cite{IoossAdel,HaragusIooss}. For application to vibratory systems, the first approaches to solve this dependency have been to compute a numerical manifold for each forcing frequency as in~\cite{JIANG2005H}, or simplify the forcing effect and consider the simplest time-dependence of the manifold, without deformations, as {\it e.g.} in~\cite{TOUZE:JSV:2006}. Using the parametrisation method and following for example~\cite{haro2006parameterization2}, a small forcing assumption allows one to limit the forcing effect to its first $\varepsilon$-order, and split the contributions of the autonomous and non-autonomous terms in separate calculations, a method that has been first applied to vibratory systems in~\cite{BreunungHaller18,HallerIsola}. In conjunction with this first $\varepsilon$-order assumption, the quantities related to the non-autonomous terms can also be solved with arbitrary order expansions. In the context of direct calculations, whereas a zero-order development for these terms is considered in~\cite{JAIN2021How}, arbitrary expansions are consistently derived in~\cite{opreni22high}, with the main consequence of an improved gain in accuracy, together with the ability to treat the case of parametric resonances for example.

The results shown in~\cite{opreni22high} on academic and engineering structures highlight the need to use a high-order expansion for the $\varepsilon$-order forcing term, bringing very important corrections to the predictions in terms of amplitudes. Nevertheless, the assumption of small forcing, which is then set to an $\varepsilon$ amplitude, limits the generality of the development. Even if the forcing amplitudes have been quantified in~\cite{opreni22high}, showing that this development allows to deal with a comfortable range of vibrations covering a variety of engineering situations, the method still has a limitation and for instance is not able to tackle superharmonic resonances. To understand this, one has to inspect the resonance relationships that arise naturally in the solution procedure, making appear the nonlinear resonances that are intimately connected to the normal form theory~\cite{IoossAdel,gucken83,Wiggins,Poincare,Dulac1912,NayfehNF}. In the  $\varepsilon^1$-order assumption used in~\cite{opreni22high,BreunungHaller18}, 
the forcing terms appear with power one in the reduced dynamics expressed with normal coordinates. The resonance relationships for this case have been  analysed in~\cite{opreni22high}, showing that, with a single master mode assumption, they reduce to the simple expression $\Omega = \gamma \omega_m$, with $\Omega$ the forcing frequency, $\gamma \geq 1$ an integer, and $\omega_m$ the eigenfrequency of the master mode. Consequently, the $\varepsilon^1$-order assumption on the forcing allows one to treat the case of primary resonance ($\gamma=1$), parametric resonance ($\gamma=2$), and 1:3 subharmonic resonance ($\gamma=3$). But to deal with a superharmonic resonance, one needs multiple occurrences of the forcing frequency in the resonance relationships, which comes in the calculation only by considering higher $\varepsilon$ orders in the development of the forcing. 

Previous developments as those reported for example in~\cite{ShawShock,JIANG2005H}, already proposed to deal with time-dependent invariant manifold, without relying on an $\varepsilon^1$-order truncation for the forcing. While in~\cite{ShawShock}, the forcing was added as an added oscillator; in~\cite{JIANG2005H} only the excitation frequency was added as a new variable to make the system autonomous. For the solution, these papers considered only the graph style parametrisation, following their previous works. While this procedure could have been considered again and enlarged to the case of the parametrisation method for invariant manifold, it appears to face two limitations. First, existence and uniqueness of the searched time-dependent invariant manifolds are not proven theoretically by mathematical theorems in the case of an eigenvalue with vanishing real part, which is the classical case under study in vibration theory to deal with harmonic forcing. Second, The initial condition of the auxiliary variable representing forcing is unitary, see {\em e.g.}~\cite{ShawShock}, which can be problematic in the context of solutions based on asymptotic expansions. Our goal here is thus to stick to the case where existence and uniqueness have already been proven, which is only possible by treating the forcing as a small parameter and computing perturbation of the autonomous invariant manifold, as shown for instance in~\cite{haro2006parameterization2,Haller2016}, with the added value that the expansions will not  be stopped to $\varepsilon^1$ order for the non-autonomous term. For the sake of completeness, the comparison between these two solution strategies is further developed in~\ref{app:nonautreat}.

The aim of this article is thus to propose an arbitrary order expansion for the direct parametrisation method of invariant manifolds. 
To reach this goal, a new version of the algorithm is proposed, in which the $\varepsilon$ development of the forcing are pushed to arbitrary order by treating the non-autonomous term as an additional parametrising variable. This embedding of the non-autonomous term into the original framework of the autonomous algorithm allows for little modifications to the structure of the algorithm itself. At the same time, the range of applicability is significantly extended, thus providing a more general method, applicable to FE problems with higher forcing levels than the current state of the art. All the developments presented herein have been implemented in a new version of the code \texttt{MORFE}, thus extending the range of the previous versions~\cite{vizza21high,opreni22high,morfe}.

\section{Direct parametrisation for non-autonomous problems}
\label{sec:methods}

In this section, a reduction method proposing a direct parametrisation for invariant manifolds of non-autonomous systems is detailed.
First, a preamble is devoted to presenting the proposed strategy to cope with the non-autonomous term, and how it compares with the current state of the art. Then the case of first-order dynamical systems is developed. 
Finally, the application to second-order nonlinear vibratory systems is detailed, where the features of the formulation are used to decrease the computational burden.

\subsection{Treatment of the non-autonomous terms in the context of the parametrisation method}
\label{subsec:preamble}

The aim of this section is to present a methodology that can be used to deal with the non-autonomous terms for a forced dynamical system. In order to make the presentation light, let us simplify the framework by considering a $D$-dimensional first-order dynamical systems with a single external forcing term:
\begin{equation}
\label{eq:FODE}
    \B \dot{\y} = \A\y +\Q(\y,\y) 
    + \varepsilon \C \e^{\tilde{\lambda} t},
\end{equation}
where the state-space vector $\y\in\mathbb{C}^D$. The real-valued matrix $\B$ is not assumed to have any general property; in particular, it might be a singular matrix such that some lines of the problem could not involve time derivatives, meaning that the development also considers the case of differential-algebraic equations (DAE), a case that has already been treated in~\cite{li2023constrained}, to which the interested reader is referred for more theoretical details.
On the other hand, the real-valued matrix $\A$ is assumed to have full rank. The nonlinearity is  given by a smooth analytical function in quadratic form $\Q(\y,\y)$. This choice not only makes the derivations simpler, but also allows, thanks to the quadratic recast~\cite{COCHELIN2009,KARKAR2013,Guillot:recast}, to treat any analytical nonlinearity by augmenting the size of the system and adding new coordinates. Note that performing quadratic recast is not always the best choice from a computational point of view and, in the case of nonlinearities that are easy to incorporate in the algorithm, such as cubic ones, it might be better to keep them on the right-hand side, as shown in Section~\ref{subsec:mechsyso2}.
The forcing term has a time dependence $\e^{\tilde{\lambda} t}$ with $\tilde{\lambda}$ the forcing value. Constant forcing term is here not considered since its effect is to change the location of the fixed point. It is hence assumed that Eq.~\eqref{eq:FODE} describes the dynamics in the vicinity of the origin, and in case a static solution needs to be considered, then a change of coordinate has to be first applied, as done for example in~\cite{opreniPiezo,Martin:rotation} in the same context. The forcing vector is represented by the vector $\C$ indicating its direction, and a small parameter $\varepsilon$. For the sake of simplicity we assumed a single forcing term, but the case of multiple ones can be treated in the same way. All the methods presented here are local in nature and rely on asymptotic expansions, motivating the introduction of this $\varepsilon$ scaling for the forcing. 

In the autonomous case, the main idea of the parametrisation method for invariant manifolds is to introduce a nonlinear mapping $\bar{\W}(\bar{\z})$ relating the original coordinates of the dynamical system $\y$ to a new {\it normal} coordinate $\bar{\z}$ describing the dynamics on the associated invariant manifold of dimension $d\ll D$. Note that $\bar{\z}$ is a vector collecting the normal coordinates of the master modes only, and is thus of dimension $d$. Together with this unknown mapping, the reduced dynamics $\bar{\f}(\bar{\z})$ governing the time evolution of $\bar{\z}$ along the embedding is also introduced. Elimination of time provides the invariance equation from which homological equations are deduced by separating orders, hence providing recursive solutions~\cite{Cabre3,Haro}.

In the presence of an external excitation, the autonomous invariant manifold becomes time-dependent, as it starts to oscillate and deform under the action of the forcing term. Provided that the forcing is small in the mathematical sense and under appropriate non-resonance conditions, a non-autonomous invariant manifold is proven to exist, see {\it e.g.}~\cite{haro2006parameterization2,Haller2016}.
This manifold can be described by a time-dependent parametrisation
$\W^{(\varepsilon)}(\bar{\z},t)$ and a time-dependent reduced dynamics $\dot{\bar{\z}} = \f^{(\varepsilon)}(\bar{\z},t)$, which satisfy the invariance equation obtained by substituting $\y$ with $\W^{(\varepsilon)}(\bar{\z},t)$ into \eqref{eq:FODE}, thus leading to
\begin{equation}
\label{eq:invariance_time}
    \B \nabla_{\bar{\z}} \W^{(\varepsilon)}\f^{(\varepsilon)} 
    +
    \B \partial_t\W^{(\varepsilon)} 
    = 
    \A \W^{(\varepsilon)}
    + 
    \Q(\W^{(\varepsilon)},\W^{(\varepsilon)})
    +
    \varepsilon \C \e^{\tilde{\lambda} t}.
\end{equation}

Although the mathematical results on existence and persistence are derived in an abstract setting~\cite{haro2006parameterization2}, meaning they apply to arbitrary order expansion in $\varepsilon$, in all the previous applications to forced systems,
see for example~\cite{haro2006parameterization2,BreunungHaller18,HallerIsola,JAIN2021How,opreni22high}, the unknown nonlinear mappings and reduced dynamics were introduced as a two terms expansion limited to the $\varepsilon^1$ order:
\begin{subequations}\label{eq:ziparamepsdev}
\begin{align}
    \y &= \bar{\W}(\bar{\z}) + \varepsilon \hat{\W}(\bar{\z},t) + \mathcal{O}(\varepsilon^2), \label{eq:ziparamepsdeva}\\
    \dot{\bar{\z}}& = \bar{\f} (\bar{\z}) + \varepsilon \hat{\f}(\bar{\z},t) + \mathcal{O}(\varepsilon^2).\label{eq:ziparamepsdevb}
\end{align}
\end{subequations}
Thanks to this expansion, the first terms $ \bar{\W}$ and $\bar{\f}$ respectively refer to the nonlinear mapping and the reduced dynamics of the autonomous system, whereas the time-dependent part due to the forcing is embedded in $\hat{\W}$ and $\hat{\f}$. The invariance equation could then be split by order, see {\it e.g.}~\cite{BreunungHaller18,opreni22high}, leading to
\begin{align}
    &
    \label{eq:invariance_eps0}
    \varepsilon^0:\qquad
    \B \nabla_{\bar{\z}} \bar{\W} \bar{\f}
    = 
    \A \bar{\W}
    + 
    \Q(\bar{\W},\bar{\W}),
    \\
    &
    \label{eq:invariance_eps1}
    \varepsilon^1:\qquad
    \B \nabla_{\bar{\z}} \bar{\W} \hat{\f}
    +
    \B \nabla_{\bar{\z}} \hat{\W} \bar{\f}
    +
    \B \partial_t\hat{\W}
    = 
    \A \hat{\W}
    + 
    \Q(\bar{\W},\hat{\W})
    + 
    \Q(\hat{\W},\bar{\W})
    +
    \varepsilon \C \e^{\tilde{\lambda} t}.
\end{align}
It is clear how the $\varepsilon$ scaling allows to separate the contributions of autonomous and non-autonomous terms, so they can be computed one after the other. The $\varepsilon^0$ problem corresponds to the autonomous system, and once it is solved for $\bar{\W}$ and $\bar{\f}$, the $\varepsilon^1$ problem can be computed to yield $\hat{\W}$ and $\hat{\f}$.
Besides, the $\varepsilon^1$ problem is a set of time-dependent linear ordinary differential equations in the unknown variables~\cite{opreni22high}, such that 
the solution to the exponential forcing can be written as an exponential, yielding
\begin{subequations}
\begin{align}
    \y &= \bar{\W}(\bar{\z}) + \varepsilon \hat{\W}(\bar{\z})\e^{\tilde{\lambda} t} + \mathcal{O}(\varepsilon^2), \\
    \dot{\bar{\z}}& = \bar{\f} (\bar{\z}) + \varepsilon \hat{\f}(\bar{\z})\e^{\tilde{\lambda} t} + \mathcal{O}(\varepsilon^2).
\end{align}
\label{eq:last_forcing_paper}
\end{subequations}

The main advantage of this processing of the forcing term, as underlined in~\cite{opreni22high}, is that the structure of the homological equations resulting from the splitting between $\varepsilon^0$ and $\varepsilon^1$ orders, is the same, such that extending the computations to non-autonomous systems is easily attainable. However, the limitation to $\varepsilon^1$-order is clear. Even though the structures of the next $\varepsilon^p$ orders also share the same structure, recursive coding of the next orders appears quite cumbersome. Numerical examples shown in~\cite{opreni22high} highlight that $\varepsilon^1$-order truncation together with arbitrary order expansion to solve for $\hat{\W}$ and $\hat{\f}$, already provides excellent improvement covering wide ranges of applications. However, for some extreme conditions, this might not be sufficient. 
Besides, the main limitation of this approach is to provide in the reduced dynamics only terms that are proportional to the forcing amplitude at power 1, 
which comes with the fact that, in the resonance relationships, the forcing frequency appears with a single occurrence and not as a summed term that can create multiples of $\tilde{\lambda}$. While this can be used advantageously in the context of parametrised ROMs with forcing, see {\it e.g.}~\cite{Martin:rotation} for an example, it is theoretically speaking a limitation since hindering the method to deal with superharmonic resonances, where higher $\varepsilon$ orders are then needed.


The aim of the present development is to bypass the first-order $\varepsilon$ truncation relative to the non-autonomous term and push the $\varepsilon$ expansion up to arbitrary order. To do so, it is convenient to introduce the following notation for the autonomous terms and the first-order terms in $\varepsilon$:
\begin{subequations}
\begin{align}
    & \W^{(\varepsilon^0)}(\bar{\z}) = \bar{\W}(\bar{\z}), \\
    & \W^{(\varepsilon^1)}(\bar{\z}) = \hat{\W}(\bar{\z}),\\
    & \f^{(\varepsilon^0)}(\bar{\z}) = \bar{\f} (\bar{\z}), \\
    & \f^{(\varepsilon^1)}(\bar{\z}) = \hat{\f}(\bar{\z}).
\end{align}
\end{subequations}

The generalisation of Eq.~\eqref{eq:last_forcing_paper} for higher order expansion in $\varepsilon$ can be written as
\begin{subequations}
\begin{align}
    \y &= \W^{(\varepsilon^0)}(\bar{\z}) + \varepsilon \W^{(\varepsilon^1)}(\bar{\z})\e^{\tilde{\lambda} t} +
    \varepsilon^2 \W^{(\varepsilon^2)}(\bar{\z})\e^{2\tilde{\lambda} t}
    + 
    \dots +
    \varepsilon^{{o}} \W^{(\varepsilon^{{o}})} (\bar{\z})\e^{{o}\tilde{\lambda} t}+
    \mathcal{O}(\varepsilon^{{o}+1}), \\
    \dot{\bar{\z}}& = \f^{(\varepsilon^0)}(\bar{\z}) + \varepsilon \f^{(\varepsilon^1)}(\bar{\z})\e^{\tilde{\lambda} t}
    + \varepsilon^2 \f^{(\varepsilon^2)}(\bar{\z})\e^{2\tilde{\lambda} t}
    + 
    \dots +
    \varepsilon^{{o}} \f^{(\varepsilon^{{o}})} (\bar{\z})\e^{{o}\tilde{\lambda} t}+
    \mathcal{O}(\varepsilon^{{o}+1}).
\end{align}
\label{eq:epsilon_expansion}
\end{subequations}
In fact, at any order ${p}$ of the $\varepsilon^{{p}}$ expansion, the dependence on time
of both $\W^{(\varepsilon)}(\bar{\z},t)$ and $\f^{(\varepsilon)}(\bar{\z},t)$, will always be conveyed by the corresponding harmonic $\e^{{p}\tilde{\lambda} t}$. This is due to the fact that the lower orders only create terms in the corresponding harmonic, and that the solution of the linear ordinary differential equations (ODE) to the exponential forcing can be written as an exponential. 

Given the sound mathematical results  on the existence and persistence of a parametrisation for the non-autonomous invariant manifold described by Eqs.~\eqref{eq:epsilon_expansion}~\cite{haro2006parameterization2,Haller2016}, the main objective of this work is to provide an algorithm to compute this parametrisation in a fully automated manner, for arbitrary order in $\varepsilon$.
To do so, the first step is to notice that every instance of $\varepsilon$ always appears together with the exponential term elevated at the same integer power. We can then write Eq.~\eqref{eq:epsilon_expansion} more compactly if we introduce the {\it dummy} variable
\begin{equation}
\label{eq:tilde_z_definition}
    \tilde{z}(t) \doteq \varepsilon\e^{\tilde{\lambda} t},
\end{equation}
so that we can write Eqs.~\eqref{eq:epsilon_expansion} as
\begin{subequations}
\begin{align}
    \y &= \sum_{{p}=0}^{{o}} \W^{(\varepsilon^{{p}})} (\bar{\z})
    \tilde{z}^{{p}}+
    \mathcal{O}(\varepsilon^{{o}+1}), \\
    \dot{\bar{\z}}& =
    \sum_{{p}=0}^{{o}} \f^{(\varepsilon^{{p}})} (\bar{\z})
    \tilde{z}^{{p}}+
    \mathcal{O}(\varepsilon^{{o}+1}),
\end{align}
\label{eq:epsilon_expansion_compact}
\end{subequations}
and the invariance equation~\eqref{eq:invariance_time} as
\begin{equation}
\label{eq:invariance_tilde_z}
    \B \nabla_{\bar{\z}} \W^{(\varepsilon)}\f^{(\varepsilon)} 
    +
    \B \partial_{\tilde{z}}\W^{(\varepsilon)} 
    \dot{\tilde{z}}
    = 
    \A \W^{(\varepsilon)}
    + 
    \Q(\W^{(\varepsilon)},\W^{(\varepsilon)})
    +
    \C \tilde{z}.
\end{equation}

The following step is to notice that $\tilde{z}$ can be algorithmically treated  in the same way the other normal coordinates in $\bar{\z}$ are. In fact, recalling that all the $\W^{(\varepsilon^{{p}})} (\bar{\z})$ and $\f^{(\varepsilon^{{p}})} (\bar{\z})$ terms will be sought as a Taylor expansion in $\bar{\z}$, it is clear from Eqs.~\eqref{eq:epsilon_expansion_compact} that treating $\tilde{z}$ as an additional normal variable is an automated way of expanding in both $\bar{\z}$ and $\varepsilon$ at the same time. Moreover, it is possible to see that the two terms on the left-hand side of the invariance equation~\eqref{eq:invariance_tilde_z} share the same structure, provided that one defines a dynamics for $\tilde{z}$. This is trivial to do as it follows from the definition of the dummy variable that $\dot{\tilde{z}}=\tilde{\lambda}\tilde{z}$.

To incorporate $\tilde{z}$ as an additional normal variable, let us define the augmented vector of normal coordinates as
\begin{equation}\label{eq:defznewgood}
    \z = 
    \begin{bmatrix}
        \bar{\z}\\
        \tilde{z}
    \end{bmatrix},
\end{equation}
We can then look for a parametrisation of the non-autonomous invariant manifold described by the nonlinear mapping
\begin{equation}\label{eq:ziNLmappinggoodone}
    \y = \W(\z),
\end{equation}
relating the physical degrees of freedom $\y$ to the augmented normal coordinate $\z$. At the same time, we define the reduced dynamics on the manifold as
\begin{equation}\label{eq:zireduceddyngood}
    \dot{\z} =  \f(\z).
\end{equation}
The unknown functions $\W$ and $\f$ will be searched for as polynomials in the augmented normal coordinates $\z$, and they will have to satisfy, up to a desired order, the invariance equation that in these settings read
\begin{equation}
\label{eq:invariance_compact}
    \B \nabla_{\z} \W(\z)\f(\z) 
    = 
    \A \W(\z)
    + 
    \Q(\W(\z),\W(\z))
    +
    \C \tilde{z},
\end{equation}
where the time-dependency is embedded into the dummy variable $\tilde{z}$. 

Comparing this invariance equation to the autonomous one, one can see that they share the same structure, with the exception of two key differences: (i) an additional term $\C\tilde{z}$ appears on the right hand side, which should be taken care of when solving the first order, as it will be discussed in Section~\ref{sec:order1}; (ii) the last entry of the vector $\f(\z)$ is not actually unknown. In fact, the dummy variable $\tilde{z}$ is only {\it treated} in the algorithm as an additional normal variable but its value is not arbitrary as the others in $\bar{\z}$. Since $\tilde{z}$ is well defined by its original definition from Eq.~\eqref{eq:tilde_z_definition}, its dynamics cannot be altered by the algorithm. 
Hence we have to impose
\begin{equation}
    f_{d+1}(\z) = \tilde{\lambda}\tilde{z},
\end{equation}
which means that we have to take care to set to zero all the other polynomial coefficients of $f_{d+1}(\z)$.

Besides taking care of these differences, from an algorithmic point of view, the same routines used for the autonomous problem can be readily extended to treat the non-autonomous one, thanks to the augmentation of the normal coordinates with the dummy variable $\tilde{z}$. Moreover, it is clear that, as compared to the present state of the art, the algorithm allows to solve, in a fully automated manner, any arbitrary order of the $\varepsilon$ expansion. 

Upon solving for these polynomial expansions, $\W(\z)$ and $\f(\z)$, the original value of the dummy variable has to  be substituted back into the equations, thus making the explicit dependence on time appear again. In particular, the parametrisation originally sought will be simply given by
\begin{align}
    &\W^{(\varepsilon)}(\bar{\z},t) 
    =
    \W(\begin{bmatrix}
        \bar{\z}\\\varepsilon e^{\tilde{\lambda}t}
    \end{bmatrix})
\\
    &\f^{(\varepsilon)}(\bar{\z},t) 
    =
    \f_{[1:d]}(\begin{bmatrix}
        \bar{\z}\\\varepsilon e^{\tilde{\lambda}t}
    \end{bmatrix})
\end{align}

As already mentioned, we restricted the presentation to the case of a single forcing term but the case of multiple ones can be handled in the same automated way by including as many dummy variables as needed. This procedure will be exemplified in Section~\ref{subsec:mechsyso2}.

As a closing remark, the computational procedure to derive time-dependent invariant manifolds could have followed the standard technique to render a non-autonomous system autonomous, an idea already used in~\cite{ShawShock,JIANG2005H}. However, as mentioned in the introduction, this approach suffers from a lack of a theoretical result ensuring the existence and uniqueness of the sought manifold. Nevertheless, a detailed comparison is provided in~\ref{app:nonautreat} for completeness.

\subsection{Eigenproperties of the first-order system}
\label{sec:eigendefo1}

Before addressing the application of the proposed methodology to cope with a non-autonomous system in the framework of the parametrisation method for invariant manifolds, the linear properties of the first-order system, are first recalled. Since no special properties of the matrices $\B$ and $\A$ have been assumed (except that $\A$ has rank $D$), right and left eigenvectors for the direct and adjoint problems are needed. The right eigenvectors $\Y_s$, $\forall\, s \; \in [1,D]$ are associated to the eigenvalues $\lambda_s$ such that
\begin{equation}\label{eq:righteigendef}
    \forall s = 1,\hdots,D:\quad \left( \lambda_s \B - \A  \right) \Y_s = \0.
\end{equation}
The left eigenvectors are denoted as $\X_s$ and are defined through
\begin{equation}\label{eq:lefteigendef}
    \forall s = 1,\hdots,D:\quad \X_s^{\star} \left( \lambda_s \B - \A  \right) = \0,
\end{equation}
where $\X_s^{\star} = \bar{\X}_s^t $  is the conjugate-transpose operation, also referred to as Hermitian transpose. Let us denote as $\X_{tot}$ and $\Y_{tot}$ the two matrices that gather respectively the $D$ left and right eigenvectors:
\begin{subequations}\label{eq:defmatrixofeigenvect}
\begin{align}
    \X_{tot} &= \begin{bmatrix}
        \X_1 & \X_2 & \hdots & \X_D
    \end{bmatrix}, \label{eq:defmatrixofeigenvectX}\\
    \Y_{tot} &= \begin{bmatrix}
        \Y_1 & \Y_2 & \hdots & \Y_D
    \end{bmatrix}, \label{eq:defmatrixofeigenvectY}
\end{align}
\end{subequations}
and $\Lambdavec$ the matrix containing the $D$ eigenvalues on the diagonal:
\begin{equation}\label{eq:defBiglambda}
    \Lambdavec = \mathrm{diag} (\lambda_1, \hdots, \lambda_D).
\end{equation}


In the case the $\B$ matrix is full rank, then an ODE system is at hand and it is assumed to be dissipative. The negative real parts of the eigenvalues are sorted by decreasing values, the less damped mode first, as:
\begin{equation}
    \label{eq:negative_eigs}
    \text{Re}[\lambda_D] \leq
    \text{Re}[\lambda_{D-1}] \leq
    \dots \leq
    \text{Re}[\lambda_{1}] < 0.
\end{equation}

Otherwise, if $\B$ is singular (and hence the original system is a DAE), the eigenvalues relating to the algebraic equations are infinite. In such a case, we assume that the equivalent ODE system obtained via an index reduction technique is dissipative and its eigenvalues are all negative, as proposed in~\cite{li2023constrained}. The stability of the fixed point is, in fact, a required condition of the parametrisation method~\cite{haro2006parameterization2,Haro,Haller2016}.

The two bases of eigenvectors can be normalised in a very general manner with 
\begin{subequations}\label{eq:normalizedFOS}
    \begin{align}
        \X_{tot}^{\star} \B \Y_{tot} & = \D, \label{eq:normalizedFOSa}\\ 
         \X_{tot}^{\star} \A \Y_{tot} & = \Lambdavec \D,\label{eq:normalizedFOSb}
    \end{align}
\end{subequations}
where $\D$ is a diagonal matrix with arbitrary entries that result from the choice on the  normalisation. Note that, in general, it is assumed that $\D$ is the identity matrix for the sake of simplicity. In order to keep the discussion as general as possible, this will not be assumed here. 

It is important to stress that these definitions serve the only purpose of clarifying the setting in which the method is applied, but the whole eigenvector matrix is never actually computed in the algorithm. In fact, all equations are written in physical coordinates and the only eigenvectors needed will be the ones chosen as master.

\subsection{First-order system with forcing}
\label{subsec:firstorderdynsys}

In this section, the complete derivation of the direct parametrisation method for systems written in first-order form and accounting for an external forcing term, is detailed. 
The parametrisation method of invariant manifolds aims at giving explicit expressions for the unknown nonlinear mapping, Eq.~\eqref{eq:ziNLmappinggoodone}, and the reduced dynamics \eqref{eq:zireduceddyngood}, which are  functions of the normal coordinate $\z$. In the present context of a non-autonomous system, $\z$ is of dimension $d+1$. The first $d$ coordinates correspond to the master modes, while the last entry is related to the forcing dummy variable $\tilde{z}$. The $d$ master modes are selected according to the problem under study, but typically they coincide with the $d-$slowest modes of the system. In the context of nonlinear vibrating systems, a frequency selection rule complemented with a resonance check is sufficient, see {\it e.g.}~\cite{vizza21high,opreni22high,opreniPiezo,Martin:rotation} for various examples.

The starting point is the invariance equation~\eqref{eq:invariance_compact}. Both unknowns $\W$ and $\f$ are searched for as polynomial expansions of arbitrary order~$o$, and are written as
\begin{subequations}
\begin{align}
    \W (\z) &= \sum_{p=1}^o \left[ \W(\z)  \right]_p,\\
    \f (\z) &= \sum_{p=1}^o \left[ \f(\z)  \right]_p,
\end{align}
\end{subequations}
where the notation $\P{.}$ is used to refer to an arbitrary  order-$p$ term. The summations start from order 1 since no constant term is considered. The polynomial and monomial representation is here introduced following the multi-index notation~\cite{ReedNotation}. A generic term of order $p$ is written as
\begin{equation}\label{eq:polyrep}
    \left[ \W(\z)  \right]_p = \sum_{k=1}^{m_p} \W^{(p,k)} \z^{\alphavec (p,k)},
\end{equation}

where $\alphavec (p,k)$ refers to the $k$-th monomial of order $p$, $k\in[1,m_p]$, $m_p$ being the number of monomials of order $p$ in $d+1$ coordinates, i.e.:
\begin{equation}
    m_p = \left( \begin{array}{c}
        p+d\\
        p
    \end{array}\right) = \frac{(p+d)!}{p!\;d!}.
\end{equation}
A given order-$p$ monomial is thus represented by the vector $\alphavec (p,k) = \left\{ \alpha_1 \quad \alpha_2 \quad \ldots \quad \alpha_{d+1} \right\}$ of size $d+1$, where each $\alpha_j$ is such that $0\leq \alpha_j\leq p$, and collects the power associated to $z_j$, such that $\sum_{j=1}^{d+1} \alpha_j = p$. The monomial associated to $\alphavec (p,k)$ simply reads $\z^{\alphavec (p,k)} = z_1^{\alpha_1}z_2^{\alpha_2} \hdots z_{d+1}^{\alpha_{d+1}}$. 
As an example, if $\z=[z_1 \,\, z_2 \,\, z_3 \,\, z_4]^{\text{T}}$ has four entries, then $\alphavec (7,k) = \left\{ 2 \,\, 1 \,\, 3 \,\, 1 \right\}$ represents the monomial $z_1^2 z_2 z_3^3 z_4$ for a given $k$ which depends on the ordering adopted. In Eq.~\eqref{eq:polyrep}, $\W^{(p,k)}$ stands for the $D$-dimensional vector of coefficients associated to the monomial $\z^{\alphavec (p,k)}$.

The solution to the invariance equation is found by writing it at each order, giving rise to the so-called homological equation of order $p$ that can be simply written from~\eqref{eq:invariance_compact} as:
\begin{equation}\label{eq:homologicpo1g}
    \B \left[ \nabla_{\z} \W (\z) \f (\z) \right]_p = \A \left[\W (\z) \right]_p + \left[ \Q (\W (\z),\W (\z)) \right]_p + \left[\C \tilde{z}\right]_p.
\end{equation}
The order-$1$ homological equation contains the forcing and is solved first, then the arbitrary order $p$ is considered.

\subsubsection{Order-1 homological equation}
\label{sec:order1}

At order 1, Eq.~\eqref{eq:homologicpo1g} reads
\begin{equation}\label{eq:homologic1o1}
    \B \left[ \nabla_{\z} \W (\z) \f (\z) \right]_1 = \A \left[\W (\z) \right]_1  + \C \tilde{z} . 
\end{equation}
Given the specific status of the last entry of the vector $\z$, it is important to follow here how the calculation proceeds with this last term. To this purpose, let us decompose the first-order terms appearing in Eq.~\eqref{eq:homologic1o1}. Since only linear terms are selected, one can write for the two unknowns:
\begin{subequations}
    \begin{align}
        \left[  \W (\z)  \right]_1 & = \W^{(1)} \z, \\
        \left[  \f (\z)  \right]_1 & = \f^{(1)} \z .
    \end{align}
\end{subequations}
The $D\times(d+1)$ matrix of unknown coefficients $\W^{(1)}$ can be written columnwise as
\begin{equation}
    \W^{(1)} = \begin{bmatrix}
        \W^{(1,1)} & \W^{(1,2)} & \hdots & \W^{(1,d+1)} 
    \end{bmatrix} = \begin{bmatrix}
        \bar{\W}^{(1)} & \W^{(1,d+1)} 
    \end{bmatrix}.
\end{equation}
In this expression, the last column $ \W^{(1,d+1)}$ has been isolated while the first $d$ columns are set together in the $D\times d$ matrix  $\bar{\W}^{(1)}$, as it coincides with the autonomous parametrisation matrix. For the matrix of coefficients $\f^{(1)}$, one can take advantage of the fact that the last line is zero except the $(d+1,d+1)$ term which is by assumption equal to $\tilde{\lambda}$. Hence the $(d+1)\times(d+1)$ matrix of unknown coefficients $\f^{(1)}$ can be expressed as
\begin{equation}\label{eq:f1dev00}
    \f^{(1)} = \begin{bmatrix}
        \bar{\f}^{(1)} & \f^{(1,d+1)} \\
        \0                   & \tilde{\lambda}
    \end{bmatrix},
\end{equation}
where $\bar{\f}^{(1)} $ is a $d\times d$ matrix of unknown coefficients and $\f^{(1,d+1)}$ a $d\times 1$ vector. Eq.~\eqref{eq:homologic1o1} can be rewritten columnwise as
\begin{equation}\label{eq:homo111}
    \B \begin{bmatrix}
        \bar{\W}^{(1)}  \bar{\f}^{(1)} \; \; &\;  \; \tilde{\lambda} \W^{(1,d+1)} + \bar{\W}^{(1)} \f^{(1,d+1)}
    \end{bmatrix} \z = \A \begin{bmatrix}
        \bar{\W}^{(1)} \;\;  & \;\;  \W^{(1,d+1)}
    \end{bmatrix} \z  + \begin{bmatrix}
        \0 & \C
    \end{bmatrix}\z,
\end{equation}
where the last column has been isolated in each matrix. Using the definition of $\z=[\bar{\z}^t \quad \tilde{z}]^t$ given in Eq.~\eqref{eq:defznewgood}, one can also rewrite Eq.~\eqref{eq:homo111} as
\begin{equation}
\label{eq:homo111b}
    \B\bar{\W}^{(1)}\bar{\f}^{(1)}\bar{\z}
    +\B(\tilde{\lambda} \W^{(1,d+1)} + \bar{\W}^{(1)} \f^{(1,d+1)})\tilde{z} = \A\bar{\W}^{(1)}\bar{\z}+
    \A\W^{(1,d+1)}\tilde{z}+ \C\tilde{z},
\end{equation}
which clearly makes appear the part that is linked to the autonomous normal variables,
and the terms that are related to the added coordinate $\tilde{z}$. 
This last equation is true for any $(\bar{\z}, \,\tilde{z})$, hence it can be split into two different problems. Considering the first $d$ columns leads to:
\begin{equation}
    \B \bar{\W}^{(1)}\bar{\f}^{(1)} = \A \bar{\W}^{(1)}.
\end{equation}
Here one recognises the linear eigenvalue problem. To enforce tangency to the linear eigenmodes, the solution is selected as the right eigenvectors $\Y_k\in\mathbb{C}^D$ associated to the master eigenvalues $\{ \lambda_k\}_{k\in[1,d]}$, see Eq.~\eqref{eq:righteigendef}, 
such that finally:

\begin{equation}\label{eq:Wsolorder1Y}
    \forall \, k=1,\hdots,d,\quad  \W^{(1,k)} = \Y_k.
\end{equation}
Collecting the $d$ right eigenvectors in the matrix of master modes $\Y$ as
\begin{equation}\label{eq:YrightMASTER}
    \Y = \begin{bmatrix}
        \Y_1 & \Y_2 & \hdots & \Y_d 
    \end{bmatrix},
\end{equation}
one simply has $\bar{\W}^{(1)} = \Y$. Finally for the reduced dynamics linear term, one obtains:
\begin{equation}\label{eq:f1dev01}
    \bar{\f}^{(1)} = \mathrm{diag} (\lambda_1,\hdots,\lambda_d).
\end{equation}
To conclude the calculation, the last column of Eq.~\eqref{eq:homo111} leads to the following problem:
\begin{equation}\label{eq:order1supcoldp1}
    \left( \tilde{\lambda} \B - \A   \right) \W^{(1,d+1)} = \C - \B \Y \f^{(1,d+1)}.
\end{equation}
The solution to this last equation is made difficult by the fact that there are two unknowns, namely $\W^{(1,d+1)}$ and $\f^{(1,d+1)}$ for a single equation, as well as by the fact that the left-hand side term can become singular in case of a primary resonance, if the forcing value $\tilde{\lambda}$, is aligned with one master eigenvalue. This problem is however classical in the context of the parametrisation method for invariant manifold, generalises to arbitrary order, and is solved by using the different styles of solution, see {\it e.g.}~\cite{Haro} for general discussions and~\cite{vizza21high,opreni22high} for solutions operating from the physical space.

Here we will only explain how to solve Eq.~\eqref{eq:order1supcoldp1} in a direct manner, without projection to the modal space. However, for the sake of completeness, the solutions operated in the modal space are made explicit in~\ref{app:projmodalhomo}. To set apart the case where the singularity might appear, one needs to introduce $\mR^{(1,d+1)}$ as the set of modes that are in primary resonance with the forcing value $\tilde{\lambda}$. In short, $\mR^{(1,d+1)}$ contains any $r$ mode $\Y_r$ such that $\lambda_r \simeq \tilde{\lambda}$. Note that the cardinality of $\mR^{(1,d+1)}$ can be larger than 1, for the case of a degenerate eigenvalue $\lambda_r$ with multiplicity larger than 1. For all $r\in \mR^{(1,d+1)}$, the matrix $\tilde{\lambda} \B - \A $ is nearly singular, and its kernel has the dimension of the cardinality of $\mR^{(1,d+1)}$.

Following the general discussion on the choice to make in such case (see~\ref{app:projmodalhomo} for details and proof), the components of $\W^{(1,d+1)}$ that are parallel to the kernel subspace cannot be derived from Eq.~\eqref{eq:order1supcoldp1}, and they must be set to zero, which generates an additional set of equations that have to be appended to Eq.~\eqref{eq:order1supcoldp1} in order to make it solvable while imposing this vanishing condition:
\begin{equation}\label{eq:kernelperp}
    \forall \, r \in \mR^{(1,d+1)}, \quad \X_r^{\star} \B \W^{(1,d+1)} = 0.
\end{equation}
For the other non-resonant components, $r \notin \mR^{(1,d+1)}$, in this case the choice retained (see~\ref{app:projmodalhomo} for details) is to set
\begin{equation}\label{eq:vanishfnonres}
    \forall \, r \in \mR^{(1,d+1)}, \quad f_r^{(1,d+1)} = 0.
\end{equation}
The augmented solvable system combining Eqs.~\eqref{eq:order1supcoldp1}-\eqref{eq:kernelperp}-\eqref{eq:vanishfnonres} finally reads:

\begin{equation}\label{eq:augmentedorder1final}
    \begin{bmatrix}
    \tilde{\lambda}\B- \A &
    \B \Y_{\mR} &
    \0
    \\
    \X_{\mR}^\star\B &  \0 & \0\\
    \0  & \0 & \I
    \end{bmatrix}
    \begin{bmatrix}
    \W^{(1,d+1)}\\
    \f_{\mR}^{(1,d+1)}\\
    \f_{\cancel{\mR}}^{(1,d+1)}
    \end{bmatrix}
    =
    \begin{bmatrix}
    \C\\\0\\\0
    \end{bmatrix}.
\end{equation}
where the two matrices $\Y_{\mR}$ and $\X_{\mR}$ containing respectively the resonant left and right eigenvectors have been introduced to simplify notations as
\begin{subequations}\label{eq:resonantleftright}
\begin{align}
    \Y_{\mR} = \begin{bmatrix}
        \Y_{r_1} & \hdots & \Y_{r_p}
    \end{bmatrix}, \quad \forall \, r_j \, \in \, \mR^{(1,d+1)}, \\
    \X_{\mR} = \begin{bmatrix}
        \X_{r_1} & \hdots & \X_{r_p}
    \end{bmatrix}, \quad \forall \, r_j \, \in \, \mR^{(1,d+1)} ,    
\end{align}
\end{subequations}
and the notation for $\mR^{(1,d+1)}$ has been abbreviated to  $\mR$ in the subscripts.
Besides, the vector of unknown coefficients of the reduced dynamics has been split into two parts, by separating the resonant ones collecting all the indices belonging to $\mR^{(1,d+1)}$, and denoted as $\f^{(1,d+1)}_{\mR}$, from the non-resonant ones collecting the indices that do not belong to $\mR^{(1,d+1)}$, and denoted as $\f^{(1,d+1)}_{\cancel{\mR}}$.

Notice that, since the check for resonances is only done between $\tilde{\lambda}$ and the master eigenvalues $\{ \lambda_k\}_{k\in[1,d]}$, the last entry of $\f^{(1,d+1)}$ related to the dynamics of the added variable $\tilde{z}$, is included neither in $\mR$, nor in $\cancel{\mR}$. In fact, as mentioned in Section~\ref{subsec:preamble}, the value of $\f_{d+1}$ is known and cannot be altered by the algorithm. Note that this consideration holds for any order.

\subsubsection{Order-p homological equation}
\label{sec:orderp}

To conclude this section, the solution to the homological equation of order $p$ is derived. Selecting order $p$ from the invariance equation, one has:
\begin{equation}\label{eq:homologicpo1}
    \B \left[ \nabla_{\z} \W (\z) \f (\z) \right]_p = \A \left[\W (\z) \right]_p  + \left[ \Q (\W,\W) \right]_p .
\end{equation}

The idea is to write this homological equation at the level of an arbitrary monomial $\z^{\alphavec (p,k)}$ defined by the vector of integers $\alphavec (p,k) = \{ \alpha_1, \; \ldots , \; \alpha_{d+1} \}$. To do so we need to isolate the unknown vectors $\W^{(p,k)}$ and $\f^{(p,k)}$ from the known ones calculated in previous instances of the iterative procedure. Each $(p,k)$-homological will consist of a left-hand side containing the unknowns and a right-hand side containing previously calculated quantities. Since the following derivation is quite heavy, the interested reader can refer to~\ref{app:pseudocode} for more details.
Here, suffice to say that each homological equation of order $p$ at the level of the arbitrary monomial $(p,k)$ has the same structure, which reads
\begin{equation}
\label{eq:zegoodhomolgicopo1:at}
    \left( \sigma^{(p,k)} \B - \A \right) \W^{(p,k)} + \sum_{s=1}^d \B \Y_s f_s^{(p,k)} = \R^{(p,k)},
\end{equation}
where $\R^{(p,k)}$ aggregates all quantities generated by lower order monomials, the full expression of which being detailed in~\ref{app:pseudocode}. The second term appearing in Eq.~\eqref{eq:zegoodhomolgicopo1:at} is $\sigma^{(p,k)}$, defined as:
\begin{equation}\label{eq:defsigmapk}
    \sigma^{(p,k)} = \sum_{s=1}^{d+1} \alpha_s \lambda_s .
\end{equation}
This term is of particular importance and one can see that it involves multiples of each eigenvalue, including the forcing value $\tilde{\lambda}$.
This property allows generating superharmonic resonances with the master eigenvalues. For instance, the case of $\sigma^{(p,k)}=3\tilde{\lambda}$ can appear in the algorithm, thus covering the case of a 3:1 superharmonic resonance. Conversely, state-of-the-art implementations reported for example in~\cite{opreni22high,JAIN2021How,li2021periodic}, can only handle direct resonances with a unitary multiplier in front of $\tilde{\lambda}$. 
 
Let us now highlight how to solve the homological equation of order $p$ for an arbitrary monomial, Eq.~\eqref{eq:zegoodhomolgicopo1:at}. This equation is underdetermined since there are still two unknowns, respectively the nonlinear mapping coefficients $\W^{(p,k)}$  and the reduced-order dynamics coefficients $\f^{(p,k)}$. Furthermore, the term in factor in front of $\W^{(p,k)}$ might become singular when a nonlinear resonance is met. Additional equations need thus to be appended to deal with underdeterminacy, and different kinds of solutions, namely the different styles of parametrisations, are introduced.

In the modal space, Eq.~\eqref{eq:zegoodhomolgicopo1:at} 
is much simpler to handle. The interested reader might find this derivation in~\ref{app:projmodalhomo}, which follows classical choices made in the derivation of the parametrisation method for invariant manifold~\cite{Haro}. Here, only the solution operated from the physical space is explained, which  is achieved by augmenting the system to make it solvable.

Let us introduce $\mR^{(p,k)}$ the resonant set, which collects all the $r$ indexes such that the nonlinear resonance relationship $\lambda_r \simeq \sigma^{(p,k)}$ is fulfilled:
\begin{equation}\label{eq:defRpk:main}
    \mR^{(p,k)} = \{ r \in [1,d] \; | \;  \lambda_r \simeq \sigma^{(p,k)} \}.
\end{equation}
Importantly, the index $r$ here covers only the master modes and does not contain the $(d+1)$-th term related to the forcing, following the choices made at the beginning to solve the non-autonomous problem. For an index $r\in \mR^{(p,k)}$,  the matrix $(\sigma^{(p,k)} \B- \A)$ is nearly singular and its kernel has the same dimension as the cardinality of $\mR^{(p,k)}$. The components of $\W^{(p,k)}$ parallel to the kernel subspace cannot be found from Eq.~\eqref{eq:zegoodhomolgicopo1:at}. For this reason, they must be set to zero, whatever the style used, see~\ref{app:projmodalhomo} for more details. This leads to considering the added equations (the derivation of which is detailed in~\ref{app:projmodalhomo}):
\begin{equation}\label{eq:kernel_orth_main}
    \X_{r}^\star \B \W^{(p,k)} = 0, \quad \forall r\in\mR^{(p,k)}.
\end{equation}
This condition imposes $\W^{(p,k)}$ to be orthogonal to the kernel of $(\sigma^{(p,k)} \B- \A)$. In a graph style parametrisation, these equations are sufficient since the choice of vanishing the coefficients of the nonlinear transform in the modal space is always retained. This is not the case in a normal form style parametrisation, since in this case, the idea is to simplify as much as possible the reduced dynamics by vanishing the coefficients for the non-resonant monomials, leading to
\begin{equation}\label{eq:vanishfnotres}
f_{r}^{(p,k)}=0, \forall r\notin\mR^{(p,k)}.
\end{equation}

Finally, in order to propose direct computations that can be done from the physical space, and using a bordering technique that augments the size of the system to avoid singularities, leads to the following problem to be solved for a given monomial with arbitrary order $p$, by grouping Eqs.~\eqref{eq:zegoodhomolgicopo1:at}, \eqref{eq:kernel_orth_main} and~\eqref{eq:vanishfnotres}:
\begin{equation}\label{eq:solopdirectcomp}
    \begin{bmatrix}
    \sigma^{(p,k)} \B- \A &   \B \Y_{\mR} & \0
    \\
    \X_{\mR}^\star\B &  \0  & \0 \\
    \0    &   \0      & \I
    \end{bmatrix}
   \begin{bmatrix}   \W^{(p,k)}\\
   \f^{(p,k)}_{\mR}\\
   \f^{(p,k)}_{\cancel{\mR}}
    \end{bmatrix}
    =
    \begin{bmatrix}
    \R^{(p,k)}\\
    \0 \\
    \0
    \end{bmatrix},
\end{equation}
where the two matrices $\Y_{\mR}$ and $\X_{\mR}$ have been introduced following Eq.~\eqref{eq:resonantleftright}, while $\f^{(p,k)}_{\mR}$ contains the coefficients $f_r^{(p,k)}$, $\forall r \in \mR^{(p,k)}$, and $\f^{(p,k)}_{\cancel{\mR}}$ contains the coefficients $f_r^{(p,k)}$, $\forall r \notin \mR^{(p,k)}$.

Note that, by projecting the first row of the system of Eqs.~\eqref{eq:solopdirectcomp} on $\X_{r}^\star$, with $r\in\mR$, an explicit expression for the reduced dynamics coefficients can be obtained as:
\begin{equation}
     f_{r}^{(p,k)} = \dfrac{\X_{r}^\star \R^{(p,k)}}{\X_{r}^\star \B \Y_{r}}, \qquad \forall r\in\mR,
\end{equation}
with $\X_{r}^\star \B \Y_{r}$ the normalisation scalar, usually equal to 1.
This expression can be used in particular to better understand how superharmonic resonances are handled by the method. Let us take the particular case of a 3:1 superharmonic resonance, such that $\sigma^{(p,k)}=3\tilde{\lambda}$. The monomial associated with this resonance relationship is $\tilde{z}^3$. If such $\sigma^{(p,k)}$ is resonant with the master eigenvalue $\lambda_r$, then, in the dynamics of the corresponding oscillator $z_r$, a forcing term $\tilde{z}^3 = e^{3\tilde{\lambda}t}$ will appear, due to the presence of a nonzero coefficient $f_r^{(p,k)}$. The term will be then responsible for exciting the oscillator at its superharmonic resonant response. 

The resonant set defined in Eq.~\eqref{eq:defRpk:main} is valid for a normal form style parametrisation. For a graph style parametrisation (see~\cite{Haro} for general discussions and~\cite{vizza21high} for specific developments for vibratory systems), one imposes $\mR_{graph}^{(p,k)}=\{1, \hdots,d\}$, consequently the last lines in Eq.~\eqref{eq:solopdirectcomp} needs not being considered, as the set $\cancel{\mR}$ is empty.

In summary, in this section, the general equations for applying the direct parametrisation of invariant manifolds to a non-autonomous problem have been derived, based on an initial assumption proposed to treat efficiently the non-autonomous term. The invariance equation has been solved at an arbitrary order and special attention has been devoted to the treatment of the non-autonomous term. 
Special attention will now be devoted to the case of mechanical systems. 
In particular, it will be underlined how one can take advantage of the fact that the initial problem is second-order in time.

\subsection{Second-order problems for mechanical vibrations}
\label{subsec:mechsyso2}

In this section, the previous developments are adapted to the case of mechanical vibratory systems featuring geometric nonlinearity. The equations of motion to be considered read~\cite{holzapfel00,ReviewROMGEOMNL}
\begin{equation}
    \M \ddot{\U} + \C\dot{\U} +\K \U + \G(\U,\U) + \Hv(\U,\U,\U) = \F(t),
\label{eq:mechstart00}
\end{equation}
where $\M$, $\C$, and $\K$ stand respectively for the mass, damping, and stiffness matrices. $\U$ stands for the displacement vector and is assumed to be of size $N$. In the context of geometric nonlinearity, and assuming that Eq.~\eqref{eq:mechstart00} is obtained from a finite element discretisation using three-dimensional elements implementing linear elasticity with a full Lagrangian formulation described by the Green-Lagrange
strain measure, conjugated with the second Piola-Kirchhoff stress measure, then the nonlinearity is polynomial of order two and three and is exactly represented by the terms $\G(\U,\U)$ and $\Hv(\U,\U,\U)$. The forcing term at the right-hand side is assumed to represent a  harmonic forcing, with a single excitation frequency $\Omega$, and reads:
\begin{equation}
\label{eq:SODEF}
    \F (t) = \varepsilon \E_+ \e^{+i \Omega t} + \varepsilon \E_- \e^{-i \Omega t},
\end{equation}
where $\E_+$ and $\E_-$ represents the spatial distribution of the harmonics of the forcing $\pm i \Omega$.

The aim of this section is to adapt the previous developments to Eq.~\eqref{eq:mechstart00}. In particular, since the two harmonics $\pm i\Omega$ are considered, the method to make the system autonomous will consider two added coordinates. Also, a particular attention will be paid to take advantage of the fact that the initial problem is second-order in time. As done for instance in~\cite{vizza21high,opreni22high}, this property can be used to halve the size of the problems to solve. In the same spirit, we choose in this development to keep the cubic nonlinearity as it is, rather than performing quadratic recast on it.

In order to use the results from the previous section, Eqs.~\eqref{eq:mechstart00} is  rewritten as a first-order system as
\begin{subequations}\label{eq:mechsyststartorder1}
    \begin{align}
         &\M \dot{\V} + \C\V +\K \U + \G(\U,\U) + \Hv(\U,\U,\U) = \E_+ z_+ + \E_- z_-, \label{eq:mechsyststartorder1a}\\
         & \M \dot{\U} = \M \V, \label{eq:mechsyststartorder1b}\\
         & \dot{z}_+ = i \Omega z_+ ,\label{eq:mechsyststartorder1c}\\
         & \dot{z}_- = -i \Omega z_- ,\label{eq:mechsyststartorder1d}
    \end{align}
\end{subequations}
which is complemented with initial conditions on the two added coordinates representing the forcing so as to ensure that they are small and can be used in asymptotic developments: $z_+ (0) = z_-(0)=\varepsilon$. To make explicit the link with the previous section, the matrices $\B$ and $\A$ are introduced as
\begin{equation}\label{eq:mechsysdefAB}
    \B = \begin{bmatrix}
        \M & \0 \\
        \0 & \M
    \end{bmatrix}, \quad \A=\begin{bmatrix}
        -\C & -\K \\
        \M & \0
    \end{bmatrix}.
\end{equation}

The fact that the system is second-order in time finds back in Eq.~\eqref{eq:mechsyststartorder1b}, which gives a direct relationship between displacement and velocity. In order to keep the size $N$ for the problems to solve in the parametrisation method, and not extend to $2N$, this property will be fully exploited, and all the vectors will be split into two parts by separating the contributions relative to the displacement and to the velocity. This can be first done for example for the right and left eigenvectors, defined in Section~\ref{sec:eigendefo1}. Let us introduce for the mechanical problems the following notation, $\forall k \in [1,2N]$:
\begin{equation}\label{eq:splitUVonXY}
    \Y_k = \begin{bmatrix}
        \Y_k^V \\
        \Y_k^U
    \end{bmatrix}, \quad \X_k = \begin{bmatrix}
        \X_k^V \\
        \X_k^U
    \end{bmatrix}.
\end{equation}
Note that at this level of the development, no special properties about the linear modes of the mechanical problems are assumed in order to keep the presentation general and underline how the method can handle different cases.

From the definition of the right eigenvectors $\Y_k$ given in Eq.~\eqref{eq:righteigendef}, and using Eqs.~\eqref{eq:mechsysdefAB}-\eqref{eq:splitUVonXY}, one can show that the two following relationships are fulfilled
\begin{subequations}\label{eq:righteigenmecdef}
    \begin{align}
        &\Y_k^V = \lambda_k \Y_k^U,\label{eq:righteigenmecdefa}\\
        &\left( \lambda_k^2 \M + \lambda_k \C + \K  \right) \Y_k^U = \0.\label{eq:righteigenmecdefb}
    \end{align}
\end{subequations}
The first equation shows that both halves of the right eigenvectors are linked through a simple relationship whatever the properties assumed on the mass, stiffness and damping matrices, as a general consequence of the link between displacement and velocity. The second equation recovers a general property that defines the eigenvalues for the mechanical problem.

Similar expressions for the left eigenvectors can be found from Eq.~\eqref{eq:lefteigendef}:
\begin{subequations}\label{eq:lefteigenmecdef}
    \begin{align}
        &{\X_k^U}^\star \M = {\X_k^V}^\star \left(\lambda_k \M + \C \right) ,\label{eq:lefteigenmecdefa}\\
        &{\X_k^V}^\star \K = - \lambda_k {\X_k^U}^\star \M.\label{eq:lefteigenmecdefb}
    \end{align}
\end{subequations}
These equations will be used to simplify some expressions in the next development.

The parametrisation method for Eqs.~\eqref{eq:mechsyststartorder1} is introduced by defining the $(2n+2)$-dimensional normal coordinate $\z$ as
\begin{equation}
    \z = \begin{bmatrix}
        \bar{\z} \\
        z_+ \\
        z_-
    \end{bmatrix},
\end{equation}
where $\bar{\z}$ is of dimension $2n$ and groups the $n$ master modes, while the two added coordinates refer to the forcing. 
The nonlinear mapping $\W$ introduced in the previous section is split into two parts. Following the notation introduced in~\cite{vizza21high,opreni22high}, it reads:
\begin{subequations}\label{eq:mappings_compact}
\begin{align}
&\U = \Psivec (\z),
\\
&
\V = \Upsvec (\z).
\end{align}
\end{subequations}
One can note in particular that on the left-hand side, only the physical displacement and velocity vectors are involved, while the reduced dynamics is a system of size $(2n+2)$ for the normal coordinates reading
\begin{equation}
    \dot{\z} = \f (\z).
\end{equation}
As underlined in the previous sections, the two last entries of $\z$ are $z_+$ and $z_-$, and it is assumed that the appended equations,
which render the system autonomous, Eqs.~\eqref{eq:mechsyststartorder1c}-\eqref{eq:mechsyststartorder1d}, are already written with the normal coordinates. As a consequence, the last two lines of $\f$, namely $f_{2n+1}(\z)$ and $f_{2n+2}(\z)$, are linear equations reproducing Eqs.~\eqref{eq:mechsyststartorder1c}-\eqref{eq:mechsyststartorder1d} and are known. Only the first $2n$ lines are nonlinear and unknown.

The general solutions obtained for the first-order problem in the previous section can now be directly applied. Using Eq.~\eqref{eq:Wsolorder1Y}, the linear terms of the mappings~\eqref{eq:mappings_compact} are obtained as, $\forall k\in[1,2n]$
\begin{subequations}
\begin{align}
    \Psivec^{(1,k)} &= \Y_k^U,\\
    \Upsvec^{(1,k)} &= \Y_k^V,
\end{align}
\end{subequations}
where $\Y_k^U$ are the displacement-related right eigenmodes associated to Eq.~\eqref{eq:righteigenmecdefb}, while $\Y_k^V$ are related to $\Y_k^U$ thanks to Eq.~\eqref{eq:righteigenmecdefa}.

The solution to the linear order needs to be completed by defining the added columns $(\Psivec^{(1,2n+1)}, \Psivec^{(1,2n+2)} )$ and $(\Upsvec^{(1,2n+1)}, \Upsvec^{(1,2n+2)} )$  corresponding to the forcing frequencies $\pm i \Omega$. To that purpose, Eqs.~\eqref{eq:order1supcoldp1} and~\eqref{eq:augmentedorder1final} needs to be rewritten twice, one for each added forcing coordinates, namely for the case $2n+1$ corresponding to $+i\Omega$, and $2n+2$ corresponding to $-i\Omega$. Let us also introduce the two resonant sets that will gather the primary resonance as
\begin{subequations}
    \begin{align}
        \mR^{(1,2n+1)} &= \mR^{(1,+)} = \{ r\in[1,2N] \; | \; \lambda_r \simeq i\Omega \},\\
        \mR^{(1,2n+2)} &= \mR^{(1,-)} = \{ r\in[1,2N] \; | \; \lambda_r \simeq -i\Omega \}.
    \end{align}
\end{subequations}
In order to alleviate notations, the superscripts $(1,2n+1)$ and $(1,2n+2)$ are replaced by $(1,+)$ and $(1,-)$ in the sequel, following the notation introduced here for 
$\mR^{(1,+)}$ and $\mR^{(1,-)}$.

Focusing only to the $(2n+1)$-th terms for the sake of brevity, and rewriting Eq.~\eqref{eq:order1supcoldp1} using $(1,+)$ notation, leads to:
\begin{subequations}\label{eq:lintermsuppmechsys}
\begin{align}
    & i\Omega \M \Upsvec^{(1,+)} + \C \Upsvec^{(1,+)}  + \K \Psivec^{(1,+)} + \sum_{r\in \mR^{(1,+)}}  \M \Y_r^V f_r^{(1,+)} = \E_+, \label{eq:lintermsuppmechsysa}\\
    &  i\Omega \M \Psivec^{(1,+)}   = \M \Upsvec^{(1,+)} - \sum_{r\in \mR^{(1,+)}}  \M \Y_r^U  f_r^{(1,+)}.\label{eq:lintermsuppmechsysb}
\end{align}
\end{subequations}
The second equation \eqref{eq:lintermsuppmechsysb} can be used to eliminate $\Upsvec^{(1,+)}$ from the calculation and halve the size of the systems to be solved thanks to:
\begin{equation}\label{eq:linkUpsPsio1}
    \Upsvec^{(1,+)} = i \Omega \Psivec^{(1,+)} + \sum_{r\in \mR^{(1,+)}}  \Y_r^U  f_r^{(1,+)}.
\end{equation}
Finally, using Eq.~\eqref{eq:righteigenmecdefa} in Eq.~\eqref{eq:lintermsuppmechsysa} in order to eliminate $\Y_r^V$ leads to a problem of size $N$ where the two unknowns are the displacement mapping term $\Psivec^{(1,+)}$ and the reduced-order dynamics coefficients $f_r^{(1,+)}$:
\begin{equation}\label{eq:mechsyso1solO}
    \left(-\Omega^2 \M + i \Omega \C + \K  \right) \Psivec^{(1,+)} + \sum_{r\in \mR^{(1,+)}}  \left[\left( (+i\Omega + \lambda_r)\M + \C  \right)\Y_r^U f_r^{(1,+)}  \right] = \E_+.
\end{equation}
This equation is interesting as it makes appear the usual linear primary resonance through the left-hand side matrix in front of $\Psivec^{(1,+)}$, which becomes singular in such case. 
The solution to Eq.~\eqref{eq:mechsyso1solO} follows the same lines of discussion as reported in the previous section and is thus not expanded further here. For a direct solution, the system needs to be augmented by bordering the matrix by the eigenvectors of its kernel, which can be directly done here by expanding the condition~\eqref{eq:kernelperp} as
\begin{equation}\label{eq:bordertemp}
    \forall r\in \mR^{(1,+)}, \quad 
    {\X_r^V}^\star \M \Upsvec^{(1,+)} +  {\X_r^U}^\star \M \Psivec^{(1,+)} = 0.   
\end{equation}
Using Eq.~\eqref{eq:linkUpsPsio1} to eliminate $\Upsvec^{(1,+)}$, and Eq.~\eqref{eq:lefteigenmecdefa} to eliminate ${\X_r^U}$, one can rewrite Eq.~\eqref{eq:bordertemp} as
\begin{equation}
    \forall r\in \mR^{(1,+)}, \quad 
    {\X_r^V}^\star \left[(+i\Omega + \lambda_r)\M + \C \right] \Psivec^{(1,+)} +  \sum_{s\in \mR^{(1,+)}}  {\X_r^V}^\star \M \Y_s^U  f_s^{(1,+)}= 0. 
\end{equation}
The two systems to be solved at order one in order to operate a direct solution from physical space can now be written down. To make the expressions more compact, the matrix of resonant right eigenvectors (displacement part) $\Y^U_{\mR}$ and left eigenvectors (velocity part)  $\X^V_{\mR}$, gathering the $r\in\mR^{(1,+)}$ terms, are introduced as
\begin{subequations}\label{eq:abbrev2np1}
    \begin{align}
        \Y^U_{\mR} &= \begin{bmatrix}
            \Y^U_{r_1} & \hdots & \Y^U_{r_p}
        \end{bmatrix}, \quad \forall r_j\in\mR^{(1,+)}, \\
        \X^V_{\mR} &= \begin{bmatrix}
            \X^V_{r_1} & \hdots & \X^U_{r_p}
        \end{bmatrix}, \quad \forall r_j\in\mR^{(1,+)}, 
    \end{align}
\end{subequations}
where the notation for $\mR^{(1,+)}$ has been abbreviated unambiguously to  $\mR$ in the subscripts. For the term $(1,2n+1)$ corresponding to $+i\Omega$, the system to be solved finally reads:

\begin{equation}
    \begin{bmatrix}
        -\Omega^2 \M + i \Omega \C + \K   &  
        \M\Y^U_{\mR}(i\Omega\I_{\mR}+\Lambdavec_{\mR})+\C\Y^U_{\mR}   & \0 \\
         (i\Omega\I_{\mR}+\Lambdavec_{\mR})\X^{V^\star}_{\mR}\M+\X^{V^\star}_{\mR}\C
         & \X^{V^\star}_{\mR} \M \Y^U_{\mR}   & \0 \\
         \0    &     \0     &   \I
    \end{bmatrix}
    \begin{bmatrix}
        \Psivec^{(1,+)} \\
        \f^{(1,+)}_{\mR} \\
        \f^{(1,+)}_{\cancel{\mR}}
    \end{bmatrix}
    = \begin{bmatrix}
        \E_+ \\
        \0 \\
        \0
    \end{bmatrix}.
\end{equation}
where $\I_{\mR}$ is an identity matrix having the cardinality of $\mR^{(1,+)}$,
and $\Lambdavec_{\mR}$ is a diagonal matrix with the same dimensions 
containing all the resonant $\lambda_r$. Finally, the coefficients of the reduced dynamics have been collected in $\f^{(1,+)}_{\mR}$ and $\f^{(1,+)}_{\cancel{\mR}}$ according to whether they belong to $\mR^{(1,+)}$ or not.

For the term $(1,2n+2)=(1,-)$ corresponding to $-i\Omega$, the same process can be conducted and one has simply to replace the $+i\Omega$ by $-i\Omega$ and $\mR^{(1,+)}$ by $\mR^{(1,-)}$, which is abbreviated to $\mR$ in the subscripts following~\eqref{eq:abbrev2np1}, so that finally one has: 
\begin{equation}
    \begin{bmatrix}
        -\Omega^2 \M - i \Omega \C + \K   &  
        \M\Y^U_{\mR}(-i\Omega\I_{\mR}+\Lambdavec_{\mR})+\C\Y^U_{\mR}   & \0 \\
         (-i\Omega\I_{\mR}+\Lambdavec_{\mR})\X^{V^\star}_{\mR}\M+\X^{V^\star}_{\mR}\C
         & \X^{V^\star}_{\mR} \M \Y^U_{\mR}   & \0 \\
         \0    &     \0     &   \I
    \end{bmatrix}
    \begin{bmatrix}
        \Psivec^{(1,-)} \\
        \f^{(1,-)}_{\mR} \\
        \f^{(1,-)}_{\cancel{\mR}}
    \end{bmatrix}
    = \begin{bmatrix}
        \E_- \\
        \0 \\
        \0
    \end{bmatrix}.
\end{equation}

As compared to the system to be solved for the first-order problem, see Eq.~\eqref{eq:augmentedorder1final}, one can observe that an additional term ${\X^V_{\mR}}^\star \M \Y^U_{\mR}$ appears in the centre of the matrix on the left-hand side. It comes from the elimination of the velocity-dependent terms and the process of halving the problems.

The same process can be repeated at arbitrary order $p$, the point being to rewrite Eq.~\eqref{eq:solopdirectcomp} to take advantage of the intrinsic properties of the second-order mechanical systems. To that purpose, the right-hand side term $\R^{(p,k)}$ needs to be also divided into two parts as:
\begin{equation}
    \R^{(p,k)} = \begin{bmatrix}
         \nuvec^{(p,k)} \\ \M \muvec^{(p,k)}
    \end{bmatrix}.
\end{equation}
Note that the matrix $\M$ can indeed be easily factorised from the lower part thanks to the definition of $\R^{(p,k)}$ given in Eq.~\eqref{eq:defRpkRHS} and the shape of the matrix $\B$ given in Eq.~\eqref{eq:mechsysdefAB}. Rewriting and expanding   Eq.~\eqref{eq:zegoodhomolgicopo1:at} then yields
\begin{subequations}\label{eq:homomechvia}
    \begin{align}
        & \left( \sigma^{(p,k)} \M + \C \right) \Upsvec^{(p,k)} + \K \Psivec^{(p,k)} + \sum_{r\in\mR^{(p,k)}} f_r^{(p,k)} \M \Y_r^V =  \nuvec^{(p,k)}, \label{eq:homomechviaa}\\
        & \sigma^{(p,k)} \M \Psivec^{(p,k)} - \M \Upsvec^{(p,k)} + \sum_{r\in\mR^{(p,k)}} f_r^{(p,k)} \M \Y_r^U = \M \muvec^{(p,k)},\label{eq:homomechviab}
    \end{align}
\end{subequations}
where the summed terms contains only resonant terms $r\in\mR^{(p,k)}$ since in the other case, when $r\notin\mR^{(p,k)}$, the corresponding reduced dynamics coefficient $f_r^{(p,k)}$ is set to zero. In order to lighten the notations again, and since one is only concerned here with the $k$-th arbitrary monomial of order $p$ unambiguously defined by the vector $\alphavec(p,k)$, the superscripts $(p,k)$ will be omitted in the next equations for $\mR^{(p,k)}$ and $\sigma^{(p,k)}$, which will thus be written simply as  $\mR$ and $\sigma$. The second equation~\eqref{eq:homomechviab} gives the important  relationship that relates the velocity mapping to the displacement mapping as
\begin{equation}\label{eq:linkPsiUpsgeneric}
    \Upsvec^{(p,k)} = \sigma \Psivec^{(p,k)} + \sum_{r\in\mR} f_r^{(p,k)}  \Y_r^U - \muvec^{(p,k)}.
\end{equation}
Replacing in Eq.~\eqref{eq:homomechviaa} and using Eq.~\eqref{eq:righteigenmecdefa} to eliminate $\Y_r^V$ leads to the homological equation at arbitrary order $p$ for the second-order mechanical systems, which reads:
\begin{equation}\label{eq:homomechsysdiv2}
    \left[ \sigma^2 \M + \sigma \C + \K  \right] \Psivec^{(p,k)} + \sum_{r\in\mR} \left[ \left(  \sigma + \lambda_r  \right)\M + \C  \right]\Y_r^U= \Xivec^{(p,k)},
\end{equation}
where the right-hand side term $\Xivec^{(p,k)}$ has been introduced to gather all the known terms as
\begin{equation}
    \Xivec^{(p,k)} = \nuvec^{(p,k)} + \left( \sigma \M +\C \right) \muvec^{(p,k)}. 
\end{equation}
Eq.~\eqref{eq:homomechsysdiv2} makes appear the nonlinear resonance relationship in the matrix in front of the mapping term $\Psivec^{(p,k)} $. The system needs to be augmented by bordering this singular matrix with the eigenvectors of its kernel, using Eq.~\eqref{eq:kernel_orth_main}, which can be expanded to
\begin{equation}\label{eq:bordertempgeneric}
    \forall r\in \mR, \quad  {\X_r^V}^\star \M \Upsvec^{(p,k)} + {\X_r^U}^\star \M \Psivec^{(p,k)} = 0.
\end{equation}
Using Eq.~\eqref{eq:linkPsiUpsgeneric} to eliminate $\Upsvec^{(p,k)}$, and Eq.~\eqref{eq:lefteigenmecdefa} to eliminate ${\X_r^U}$, one can rewrite Eq.~\eqref{eq:bordertempgeneric} as
\begin{equation}\label{eq:augmentmechsyso2final}
    \forall r\in \mR, \quad  {\X_r^V}^\star \left[ \left( \sigma + \lambda_r \right) \M + \C \right] \Psivec^{(p,k)} + \sum_{s\in\mR} f_s^{(p,k)}  {\X_r^V}^\star \M \Y_s^U =  {\X_r^V}^\star \M \muvec^{(p,k)}.
\end{equation}
The augmented system to be solved at arbitrary order for the second-order mechanical system that takes advantage of the relationship between displacement and velocity to halve the size, can be finally written, by grouping Eqs.~\eqref{eq:homomechsysdiv2} and~\eqref{eq:augmentmechsyso2final}, as:

\begin{equation}\label{eq:mecho2orderpfinal}
    \begin{bmatrix}
        \sigma^2 \M +  \sigma \C + \K   &  
        \M\Y^U_{\mR}(\sigma\I_{\mR}+\Lambdavec_{\mR})+\C\Y^U_{\mR}   & \0 \\
         (\sigma\I_{\mR}+\Lambdavec_{\mR})\X^{V^\star}_{\mR}\M+\X^{V^\star}_{\mR}\C
         & \X^{V^\star}_{\mR} \M \Y^U_{\mR}   & \0 \\
         \0    &     \0     &   \I
    \end{bmatrix}
    \begin{bmatrix}
        \Psivec^{(p,k)} \\
        \f^{(p,k)}_{\mR} \\
        \f^{(p,k)}_{\cancel{\mR}}
    \end{bmatrix}
    = \begin{bmatrix}
         \Xivec^{(p,k)}\\
        \X^{V^\star}_{\mR}  \M \muvec^{(p,k)}\\
        \0
    \end{bmatrix},
\end{equation}

where the two matrices ${\X^V_{\mR}}$ and $\Y^U_{\mR}$ gathers the left and right resonant eigenvectors as:
\begin{subequations}\label{eq:abbrev2np1orderp}
    \begin{align}
        \Y^U_{\mR} &= \begin{bmatrix}
            \Y^U_{r_1} & \hdots & \Y^U_{r_p}
        \end{bmatrix}, \quad \forall r_j\in\mR, \\
        \X^V_{\mR} &= \begin{bmatrix}
            \X^V_{r_1} & \hdots & \X^U_{r_p}
        \end{bmatrix}, \quad \forall r_j\in\mR. 
    \end{align}
\end{subequations}
In the process of halving the size of the system, two terms have appeared in the equations that are added to the augmented system. First in the centre of the matrix on the left-hand side, a term $\X^{V\star}_{\mR} \M \Y^U_{\mR} $, and second, the term on the right-hand side, $\X^{V\star}_{\mR} \M \muvec^{(p,k)}$, which makes the half-sized system very different from those obtained with a first-order formulation. One can also note that in the present derivation, very few assumptions have been made on the modes of the mechanical system. To make the link with the previous development reported in~\cite{vizza21high,opreni22high} where real normal modes were assumed, \ref{app:2ndorderreal} details how the system to be solved at arbitrary order $p$, Eq.~\eqref{eq:mecho2orderpfinal}, can be further simplified using classical assumptions for mechanical systems (real normal modes, symmetric matrices, damping matrix diagonalised by the normal modes).

\subsection{Remarks on the implementation}
\label{subsec:implementation}

The proposed formulation has been implemented in a julia package called \texttt{MORFE2.0} which can be downloaded from the github page of the \texttt{MORFE project}~\cite{morfe}. Using the multi-index notation with $\alphavec (p,k)$ vectors for processing the monomials, the code has been completely rewritten, and now gives a general formulation for non-autonomous problems.
The FEM discretisation is still based on 3D solid elements which can be 
either 15-node prisms or 27-nodes hexahedra.
The implementation follows strictly the derivation detailed in previous sections
and, actually, many definitions of the terms in the homological equations like 
\eqref{eq:Q:at}, \eqref{eq:N3:at} and \eqref{eq:N2pk:at}
have already been conveniently expressed using pseudo-code lines in the presentation given in~\ref{app:pseudocode}.

Some remarks are worth stressing.
First of all, even though for the sake of clarity
the derivations for order-1 and arbitrary order-p homological equations 
have been presented separately in Sections \ref{sec:order1} and \ref{sec:orderp}, the procedure is the same and can be unified, thus leading 
to a great simplification which is used to present the sequential calculations inside the code in a much more uniform manner.

Secondly, to keep the presentation as simple as possible,
only one forcing frequency has been introduced in Eq.~\eqref{eq:FODE}
and \eqref{eq:SODEF}, but it is actually easy to generalise to the sum of many forcing
frequencies. This leads to the presence of a larger number of non-autonomous variables
which can all a priori resonate with the autonomous ones. This added feature is thus already present in the julia package \texttt{MORFE2.0}.

Thirdly, as discussed in [27], in the usual case of $\tilde{\lambda}$ not fixed, the variation of $\tilde{\lambda}$ can be discarded, so the computation can be performed for a single $\tilde{\lambda}$ value. This assumption is made to meet the purpose of flexibility since proposing a single ROM equation that can be used for computing solution branches in the vicinity of the selected $\tilde{\lambda}$ value.

Lastly, it has been assumed everywhere in the paper 
that the maximum order of expansion~$o$ for the mappings and reduced dynamics
is actually the same as the order of $\varepsilon$ truncation.
This means that in the monomials 
$\z^{\alphavec(p,k)}$ both master mode coordinates $\bar{\z}$
and coordinates associated to the forcing (either $\tilde{z}$, $z^+$ or $z^-$)
can appear with any order ranging from $1$ to $o$.
However, this choice
considerably increases the number of homological equations to be solved
and might not be necessary for all applications.   Indeed, as will be shown in the next section, the order needed in terms of $\varepsilon$ development for the non-autonomous part is generally smaller than the polynomial expansion order. Additionally, the fact that with a $\varepsilon^1$ order development, only linear terms with respect to the forcing are present in the reduced dynamics, is an attractive feature that can be used to interpolate reduced dynamics, see {\it e.g.}~\cite{Martin:rotation}. For all these reasons, it has been thus decided to offer the user the possibility to select both the order of the asymptotic expansion $o$, as well as the order of the $\varepsilon$ truncation $o_{\varepsilon}$ as input parameters. Note that~\ref{app:truncation} gives a detailed discussion on the treatment of asymptotics by comparing the smallness of forcing and amplitudes, that appear independently in the process, but are related through the dynamics.

The notation $\mathcal{O}(z^o,\varepsilon^{o_{\varepsilon}})$ is adopted in the code and in what follows to denote the different maximum orders selected.
For instance, in the case $\mathcal{O}(z^6,\varepsilon^3)$, all the monomials in the expansion
will have maximum order 6 for both coordinates, while the non-autonomous variables will appear at most with power 3, meaning that 
we accept an order $\varepsilon^3$ truncation. This pragmatic choice can be interpreted as giving equal importance to the asymptotics related to amplitude and forcing. This point is further detailed in~\ref{app:truncation} where different truncation strategies, based on the analysis of the two asymptotics, are also illustrated and could have been used to define the truncation orders.

It is worth stressing that the case $\mathcal{O}(z^o,\varepsilon^{1})$
coincides with the formulation in~\cite{opreni22high}, hence the proposed derivation
represents a generalisation to that approach.
On the other hand, the notation $\mathcal{O}(z^o,\varepsilon^{0})$ corresponds to the case without the forcing. This choice can be used in conjunction with vanishing damping to compute the backbone curve of the associated conservative problem. It can also be used in a forced-damped case by simply appending the forcing on the right-hand side of the system to compute frequency-response curves. In this case, the manifold is assumed to stay undeformed under the action of the forcing, which is a limiting assumption already used for example in~\cite{TOUZE:JSV:2006,artDNF2020,AndreaROM,vizza21high}. A more rigorous way to include the forcing with the least possible number of added monomials is to include only linear monomials in $\tilde{z}$ (or $z^+,z^-$). In this way, one could keep higher-order monomials in the autonomous expansion $\varepsilon^0$, say $\mathcal{O}(z^o,\varepsilon^0)$, but only linear ones in the non-autonomous expansion, meaning $\mathcal{O}(z^1,\varepsilon^1)$. This approach is the one adopted in~\cite{JAIN2021How}, and has also been used in~\cite{opreni22high} to draw comparisons with higher orders developments in the $\varepsilon^1$ truncation, in order to show that the latter was needed to achieve convergence to the solution. In fact, this approach is also limited to the cases where the manifold does not deform too much under the forcing and no parametric forcing terms are present in the reduced dynamics. This case is very specific and as a matter of fact is denoted as $\mathcal{O}(z^o,\varepsilon^{0})$ by abuse of notation.

\section{Numerical results}
\label{sec:numresults}

The proposed methodology will be validated on three examples  to stress 
its enhanced capabilities for Reduced Order Model (ROM) generation. In all the examples, the ROM solutions are systematically compared to a reference 
computed by means of the Harmonic Balance Method (HBM) applied to the full-order problem (denoted as HBFEM in the following). 
The test cases  considered are quite simple from a geometrical standpoint as they concern straight beams and shallow arches. Instead, the emphasis will be put on the dynamical solutions and the ability of the ROM to recover solutions that were not computable with previous formulations of the non-autonomous part.
By resorting to non-dimensional forms of the input and output quantities, the improvements will be highlighted together with a check of the admissible range in terms of amplitudes.
It is worth stressing that geometrical dimensions and material properties 
are typical of microstructures (MEMS) that often provide experimental 
evidence of nonlinear phenomena and are recently stimulating important advancements in their numerical analysis.

\subsection{Straight beam with a 3:1 superharmonic resonance}
\label{sec:3:1}

\begin{figure}[ht]
	\centering
	\includegraphics[width = .9\linewidth]{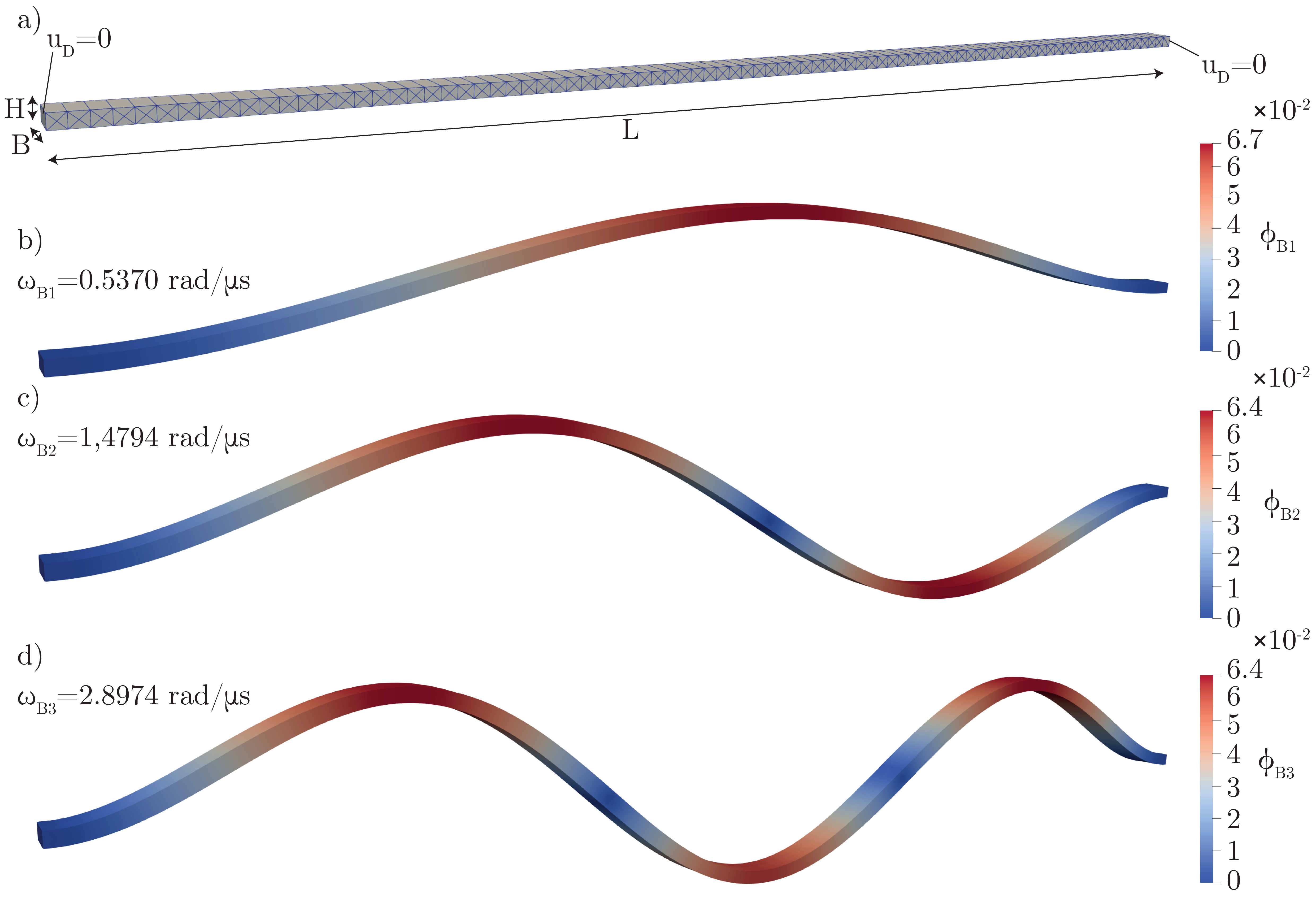}
	\caption{Clamped-clamped beam. (a) geometry $L=1000\,\mu$m, $B=24\,\mu$m, $H=10\,\mu$m. (b-d) first  three bending eigenmodes.}
	\label{fig:beam_geofig}
\end{figure}

In this section, we will address a 3:1 superharmonic resonance in 
a clamped-clamped straight beam actuated by a body load proportional to its first bending eigenshape.
The selected beam, shown in Fig.~\ref{fig:beam_geofig}(a),
is made of polycrystalline silicon, which is modeled as an
isotropic material with density $\rho=2320$\,kg/m$^3$,  
Young's modulus $E=160$\,GPa and Poisson ratio $\nu=0.22$.
The reference finite element model consists of a mesh of 
15-nodes quadratic wedge elements with 2607 nodes, yielding a system with 7821 dofs.
The eigenanalysis of the system gives the eigenmodes and the eigenfrequencies reported in Figs.~\ref{fig:beam_geofig}(b-d). Only the first three in-plane modes are reported for the sake of brevity.

The dynamics of the full-order FE model reads
\begin{equation}
    \M \ddot{\U} + \C\dot{\U} +\K \U + \G(\U,\U) + \Hv(\U,\U,\U) = \kappa \M \phivec_{B1} \cos(\Omega t),
\label{eq:mechstart}
\end{equation}
where $\phivec_{B1} $ is the mass-normalised bending mode, $\M \phivec_{B1}$ 
are nodal forces corresponding to a body force distribution proportional to the inertia 
of the mode.
In order to test the method, large values of $\kappa$ will be selected. To monitor the amplitude of the forcing with respect to the geometric nonlinearity, a non-dimensional amplitude $\epsilon$  is introduced as
\begin{equation}
\label{eq:epsload}
\epsilon=\frac{\phi}{H}\frac{\kappa}{\omega_{B1}^2}, \quad \mbox{with}
\quad \phi=\max(\|\phivec_{B1}\|).
\end{equation}
This choice is inspired by the non-dimensional form of equations in beam theory, 
see~\ref{sec:an_beam}, and follows earlier normalisation already used in~\cite{opreni22high}. Note in particular that the characteristic length selected is here the thickness $H$ in the direction of the vibration since the beam is axially constrained~\cite{ReviewROMGEOMNL}. As a general remark, "large forcing amplitude" needs to be understood here in the engineering sense. Since the parametrisation method relies on a local theory using asymptotic expansions, the non-dimensional amplitude $\epsilon$ introduced in Eq.~\eqref{eq:epsload} has to be smaller than 1, meaning that the forcing amplitudes are small in the mathematical sense. However, they are large if compared to the usual range observed in engineering literature.

To trigger the 3:1 superharmonic resonance, the excitation frequency $\Omega$ is selected in the vicinity of $\omega_{B1}/3$, where $\omega_{B1}$ is the eigenfrequency of the fundamental bending mode. The forcing  amplitude is set to $\kappa=5\,\mus/\mum^2$ which corresponds to $\epsilon=0.1170$, considering that $\phi=0.06738\mu\text{m}$. As shown in the results, this forcing will lead to vibration amplitudes in the range of 2.9$\mum$, corresponding to 
$(1/3)H$, which is sufficient to trigger nonlinear effects. 
Finally, a Rayleigh model is selected for the damping matrix $C=\alpha M+\beta K$, with $\alpha = \omega_{B1}/500, \beta =0$.

\begin{figure}[ht]
	\centering
	\includegraphics[width = .99\linewidth]{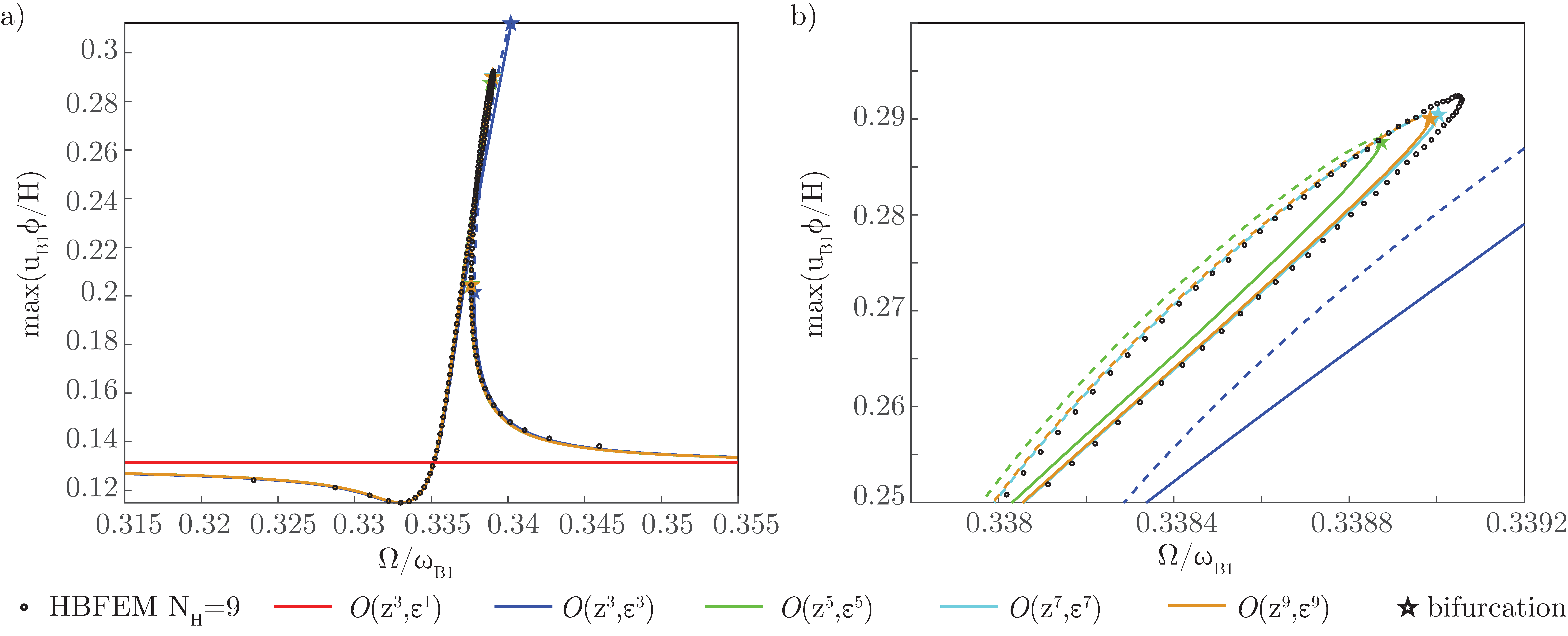}
	\caption{Frequency response curves corresponding to a 3:1 superharmonic resonance on the doubly clamped beam having a mass-proportional damping $\alpha=\omega_{B1}/500$ and an external excitation with amplitude $\kappa=5\,\mum/\mus^2$ i.e. $\epsilon=0.1170 [-] $. (a) DPIM  solution with different expansion orders compared to the HBFEM solution obtained with 9 harmonics. (b)  enlarged view of the FRC peaks. For the ROMs, the stability of the solution branches is reported with solid lines (resp. dashed lines) for stable solutions (resp. unstable solutions). The star markers pinpoint the saddle-node bifurcation points.
 }\label{fig:superharm}
\end{figure}

Fig.~\ref{fig:superharm} collects the Frequency Response Curves (FRCs) expressed in terms of the modal coordinate $u_{B1}=\phivec_{B1}^T \M \U$, normalised by the eigenmode amplitude and the characteristic length $H$.
The reference HBFEM solution has been computed with 9 harmonics together with 35 integration points for the Fourier projection. 
According to our experience, these settings provide a robust solution that is almost at convergence up to very small details. Improving the full-order solution  
would imply a computational burden difficult to manage with standard computing resources.

\begin{figure}[h!]
	\centering
	\includegraphics[width = .99\linewidth]{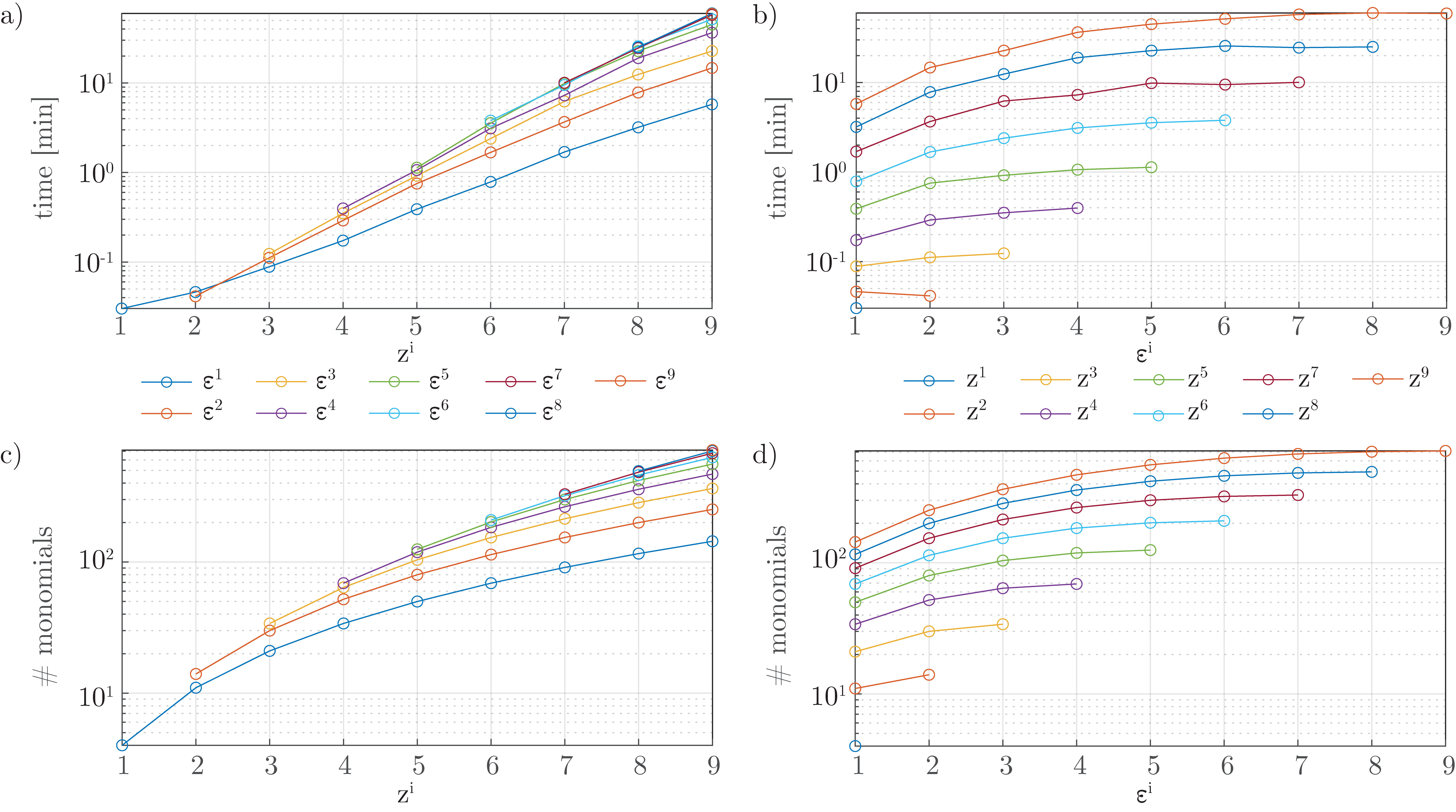}
	\caption{Computational burden to construct the ROMs in terms of time and number of monomials. (a)~Increasing orders in the parametrisation with a fixed truncation order in $\varepsilon$. (b) Increasing truncation in $\varepsilon$ order with a fixed order expansion order for the normal coordinate. (c-d) number of monomials to be computed in each case.}
\label{fig:performance}
\end{figure}

\begin{figure}[h!]
	\centering
	\includegraphics[width = .99\linewidth]{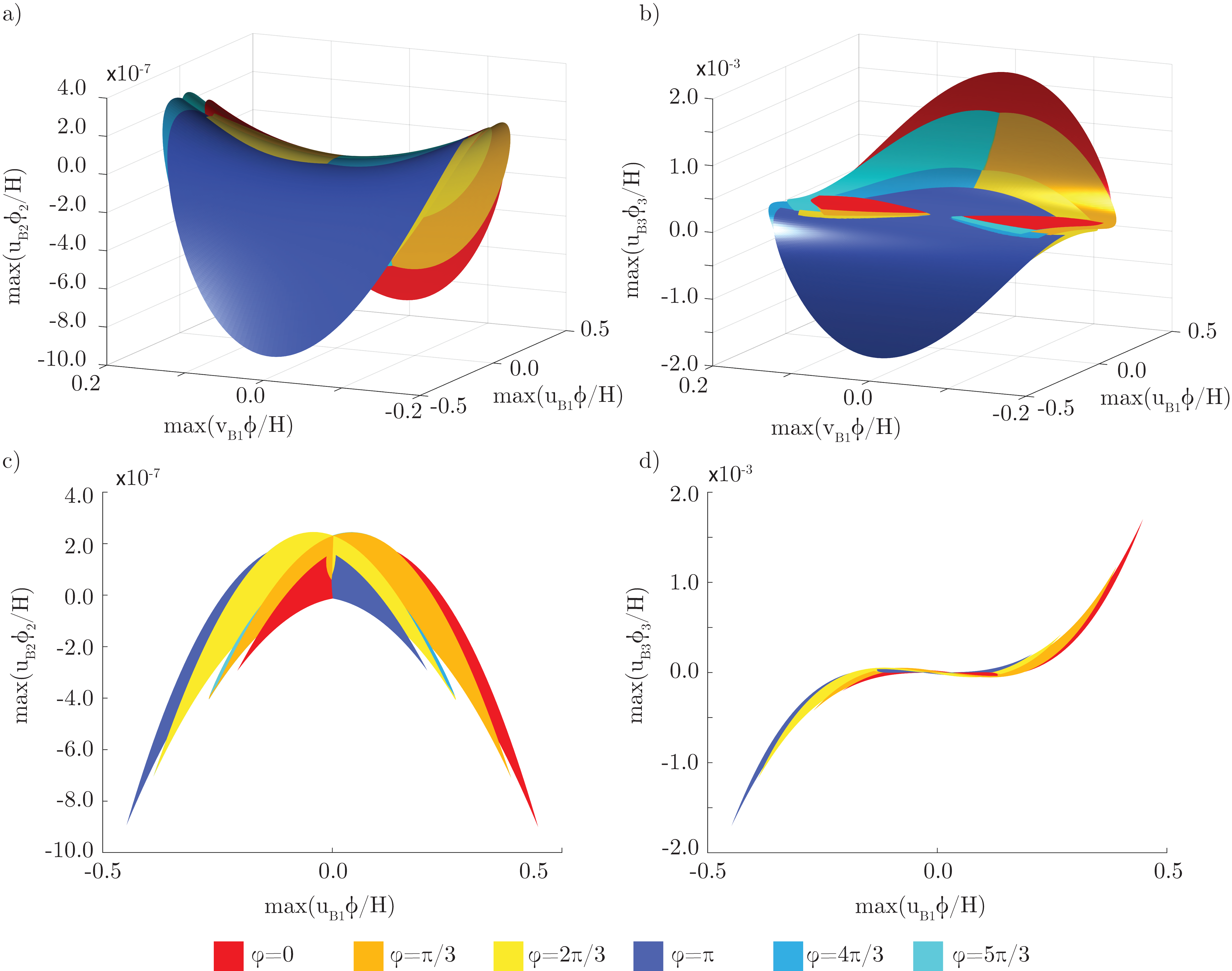}
    \caption{Graphical representation of the whisker motion with time in the case of the 3:1 superharmonic resonance in a straight beam. Single nonlinear normal mode truncation along the first bending mode (B1), parametrised with an order $\mathcal{O}(z^5,\epsilon^5)$. (a)-(c) 3-d and 2-d representation using the second bending mode as slave coordinate. (b)-(d) Using the third bending modal coordinate B3.}
 \label{fig:manifold}
\end{figure}

The ROMs are obtained by choosing a single master coordinate (corresponding to the fundamental bending mode) 
and are computed using different orders for the parametrisation
with a complex normal form style. 
The ROMs have been integrated by means of the MATLAB package MATCONT \cite{dhooge2004matcont}, which deploys continuation of periodic orbits with the collocation method. To guarantee a proper representation of the solution, 4 collocation points and 40 time intervals have been used, which results in a discretisation of the orbit with 161 nodes.

Following the introduced terminology, $\mathcal{O}(z^p,\varepsilon^q)$ defines the parametrisation and $\varepsilon$ truncation orders.  
The special case $\mathcal{O}(z^3,\varepsilon^1)$ is here used to make a direct comparison with previous formulations limited to $\varepsilon^1$ \cite{opreni22high}. As awaited, with this choice the 3:1 superharmonic resonance cannot be retrieved simply because, according to the previous theoretical developments, this resonance cannot be identified.
An $\varepsilon^3$ truncation is the minimum needed to capture the superharmonic resonance. Indeed, in Fig.~\ref{fig:superharm}  the $\mathcal{O}(z^3,\varepsilon^3)$ solution captures the resonance, but overpredicts the amplitude peak, indicating that higher orders are needed to achieve convergence. 
Finally, the cases $\mathcal{O}(z^p,\varepsilon^p)$ with $p=6$, $7$ and $9$ 
are simulated to verify the convergence of the ROM. The zoom shown in Fig.~\ref{fig:superharm}(b) highlights a slight difference between  $\mathcal{O}(z^9,\varepsilon^9)$ and the HBFEM solution. Nevertheless, the difference is very small and might also be attributed to a not fully converged HBFEM solution.

The stability of the computed branches of solutions has been reported in Fig.~\ref{fig:superharm} for the ROMs only. Indeed, computing the stability of the full-order model is computationally expensive and has not been considered here. As awaited, the stability analysis reports classical results about the solutions of the 3:1 superharmonic solution, with the presence of an unstable branch in between two saddle-node bifurcation points.

As far as the performance is considered, the $\mathcal{O}(z^9,\varepsilon^9)$ parametrisation needs less than 60 minutes to be computed, while tracing the FRCs requires only a few seconds. The $\mathcal{O}(z^5,\varepsilon^5)$ solution, which gives accurate but not fully converged results, is constructed in 1 minute. 
On the other hand, the 9 harmonics HBFEM solution requires approximately 10 minutes for each frequency value along the FRCs. As a single FRC needs 
approximately 150-200 points depending on the continuation algorithm parameters, this results in a total of 1-2 days for each curve. All the analyses have been run on a workstation with an Intel(R) Xeon(R) W-2275 CPU having 14 cores with 3.3 GHz processor base frequency, and 128 GB of RAM memory.

The computational time needed for different parametrisations is  detailed in Fig.~\ref{fig:performance} in order to give a complete overview of the construction times needed for the ROMs. The reported times correspond to the example studied, the clamped-clamped beam with 2607 nodes (7821 dofs), the reduction being performed by using a single nonlinear normal mode.  Fig~\ref{fig:performance}(a) shows that a nearly linear trend is observed when increasing the order of the asymptotic expansion for a fixed value of the $\varepsilon$ truncation for the non-autonomous part, recovering earlier results shown for example in~\cite{vizza21high}. Different orders lead to slightly different slopes. 
Fig~\ref{fig:performance}(b) underlines that, with increasing $\varepsilon$, the computational time increases at a slower pace.
This is a direct consequence of the number of terms to be computed in the parametrisation and to the truncation choices retained to define the orders $\mathcal{O}(z^p,\varepsilon^q)$, see the discussions reported in Section~\ref{subsec:implementation} and in~\ref{app:truncation}. Fig.~\ref{fig:performance}(c-d) reports the number of monomials in each truncation computed. Both figures highlight that the number of terms strongly increases with the order of the expansion in the normal coordinate $z$. On the other hand and due to the practical choice retained here to define the truncation order  $\mathcal{O}(z^p,\varepsilon^q)$,  the number of monomials to be added for each order decreases for increasing $\varepsilon$ at a fixed order for $z$.

In order to conclude this section, the motion of the parametrised invariant manifold with respect to time is displayed in Fig.~\ref{fig:manifold}, illustrating how the method automatically computes the curvatures of the reduced subspace corresponding to non-resonant coupled modes. The time-dependent invariant manifold in the case of a non-autonomous system is often referred to as a {\it whisker}~\cite{haro2006parameterization2,opreni22high}. To represent the time dependence of the whisker, different snapshots at a fixed phase are shown, which is done by fixing 6 phase values along one period of the forcing through $z^+=\sin\varphi$ and $z^-=\cos\varphi$ and $\varphi$=$0$, $\frac{\pi}{3}$, $\frac{2\pi}{3}$, $\pi$, $\frac{4\pi}{3}$, and $\frac{5\pi}{3}$.   To appreciate the curvatures of the manifold, two different slave coordinates are represented corresponding to the second and third bending modes. The plots clearly show that the most important coupling is with the third bending mode, which is a symmetric mode like the master. Whereas the deformations of the whisker are small during the oscillations, it shows a sliding motion as evidenced by the different snapshots represented in Fig.~\ref{fig:manifold}. This first observation of the whisker motion for a 3:1 resonance scenario is very different from the ones already commented at primary resonance~\cite{opreni22high}.

\subsection{Shallow arch and 2:1 superharmonic resonance}
\label{sec:2:1}

\begin{figure}[ht]
	\centering
	\includegraphics[width = .75\linewidth]{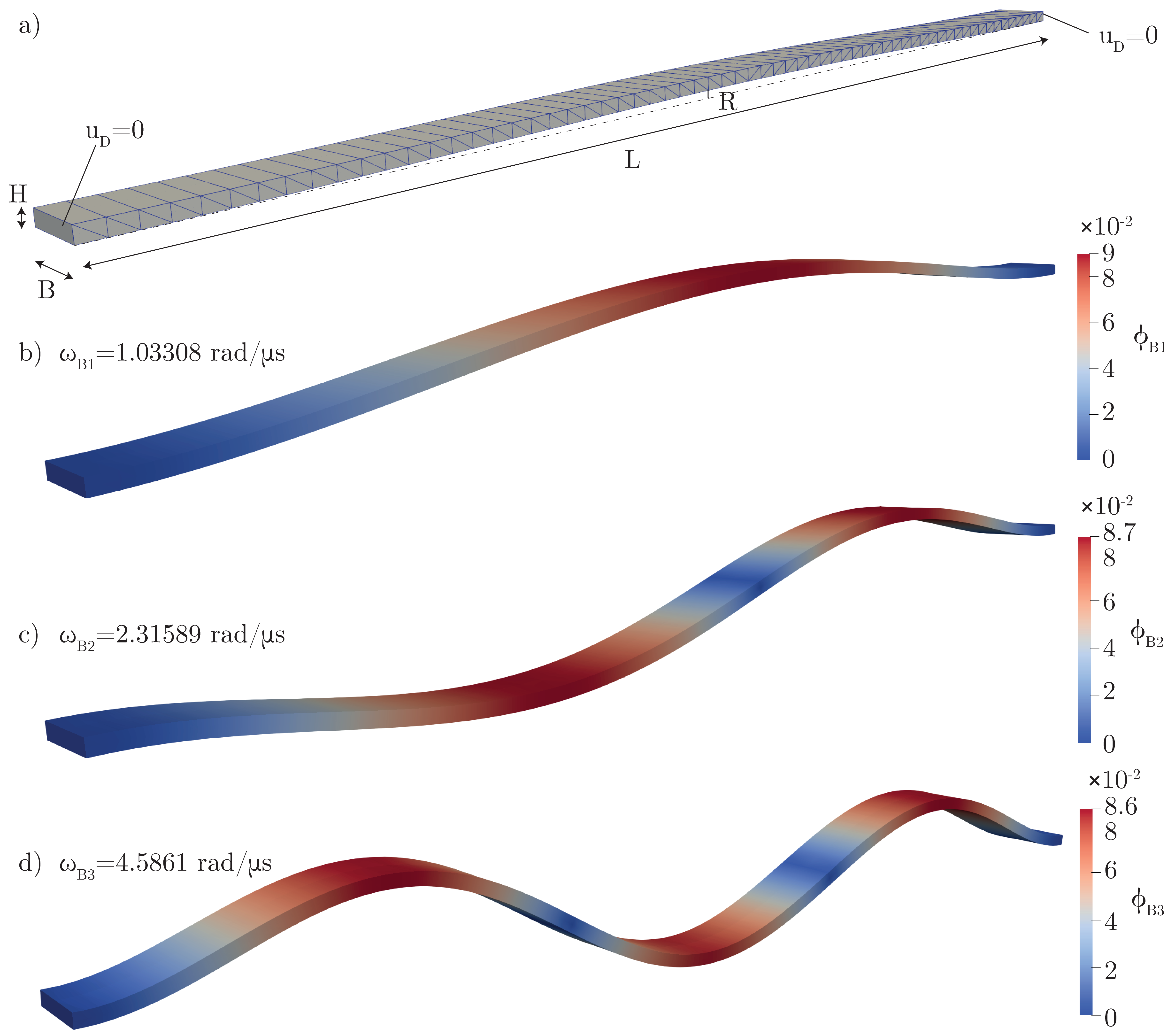}
	\caption{Sine-arched structure. (a) Geometry with length $L = 640\,\mum$, rise  
 $R = 3.84\,\mum$,  width $B=\,32\mum$ and thickness $H=\,6.4\mum$. (b-d) First three bending eigenmodes.}
	\label{fig:arch_geofig}
\end{figure}

The second selected example considers a 2:1 superharmonic resonance in a shallow arch. Whereas 3:1 superharmonics have been the subject of numerous investigations, 2:1 superharmonics excited through quadratic nonlinearity are seldom addressed. 
The structure considered is the clamped-clamped arch  shown in Fig.~\ref{fig:arch_geofig}(a) together with its dimensions. The rise is proportional to a sine function.

This beam is expected to have geometric nonlinearities with important quadratic couplings between bending modes \cite{Gobatres12,gobat2021reduced,alkharabsheh2012dynamics}.
The reference finite element model consists of a spatial discretisation made by 
15-node quadratic wedge elements with 1161 nodes, yielding a system with 3483 dofs.
Material properties are those of polycrystalline isotropic silicon with density $\rho=2320$\,kg/m$^3$,  Young's modulus $E=160$\,GPa and a Poisson ratio $\nu=0.22$. 
The first eigenmodes and eigenfrequencies are reported in Figs.~\ref{fig:arch_geofig}(b-d).
The arch is excited in the vicinity of half the eigenfrequency of the fundamental bending mode, with a body load proportional to its eigenshape.
The semi-discrete FE model reads:
\begin{equation}
    \M \ddot{\U} + \C\dot{\U} +\K \U + \G(\U,\U) + \Hv(\U,\U,\U) = \kappa_1 \M \phivec_{B1} \cos(\Omega t).
\end{equation}
with $\Omega=\omega_{B1}/2$.
As in the previous example, we adopt mass-proportional Rayleigh damping, with $\alpha=\omega_{B1}/500, \beta=0$.

\begin{figure}[ht]
	\centering
	\includegraphics[width = .99\linewidth]{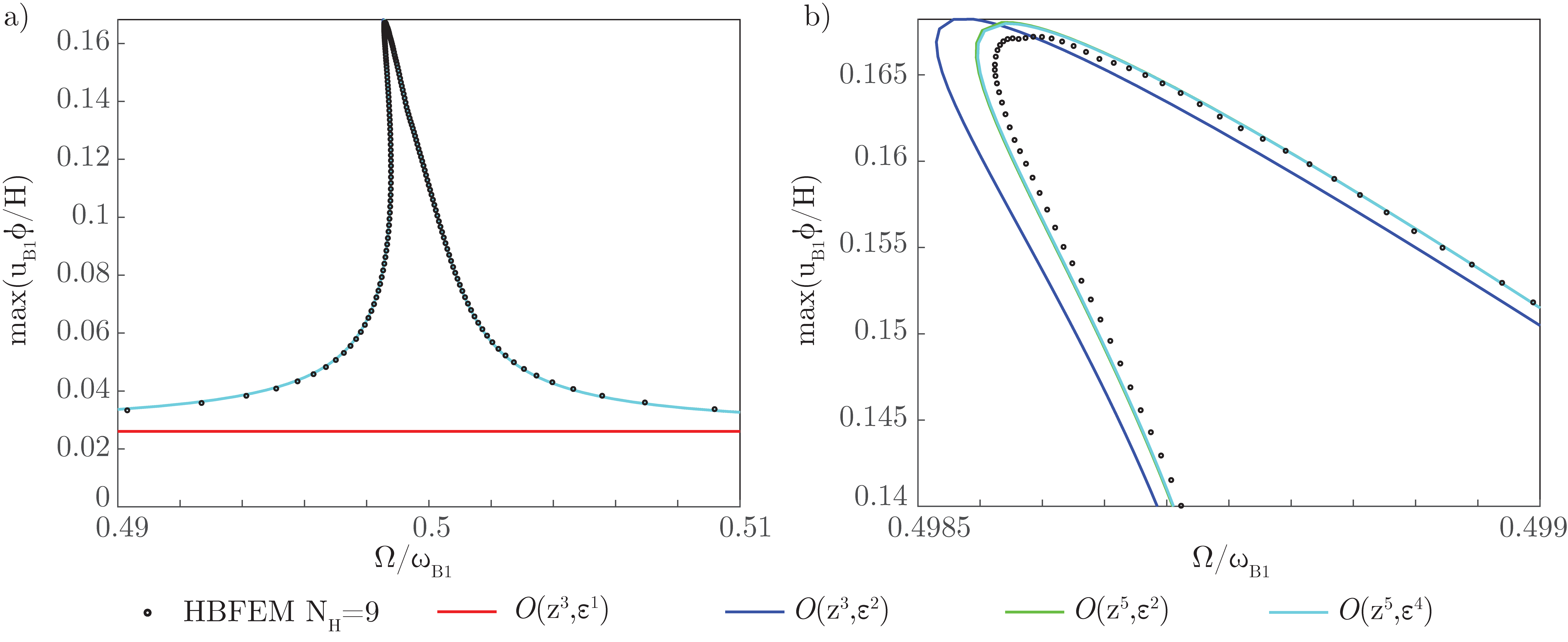}
	\caption{FRCs corresponding to a 2:1 superharmonic resonance on the arch structure having a mass proportional damping $\alpha=\omega_{B1}/500$ and an external excitation with amplitude $\kappa=1.5\,\mum/\mus^2$ i.e. $\epsilon=0.019$. (a) DPIM solution with different expansion orders compared to the HBFEM solution. (b) enlarged view of the FRCs peaks. We remark that any expansion having $\varepsilon$ strictly greater than 1 match almost perfectly the reference HBFEM.
 }	\label{fig:arch_12super}
\end{figure}

\begin{figure}[p]
	\centering
	\includegraphics[width = .99\linewidth]{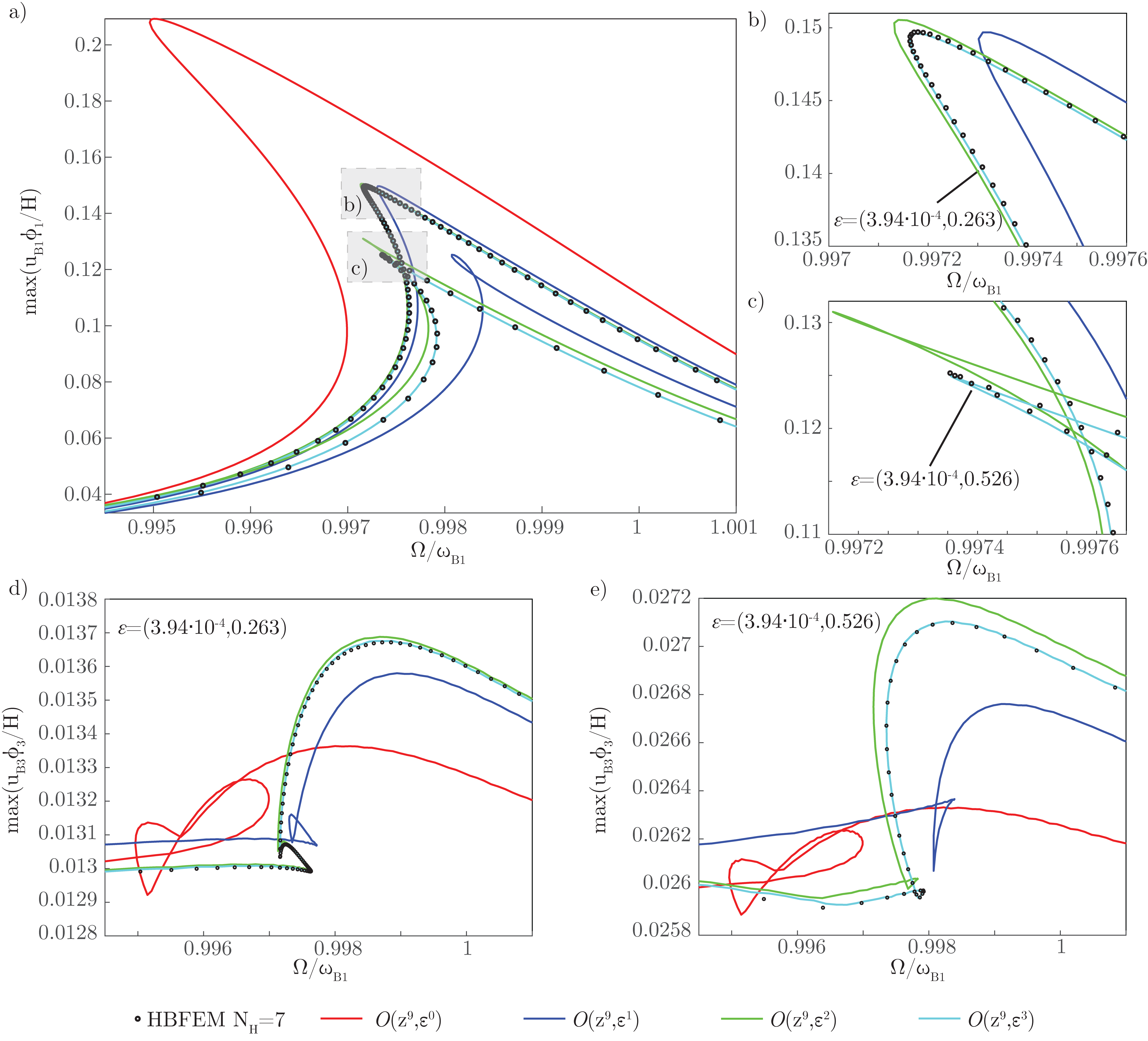}
	\caption{Multimodal forcing of a shallow arch: performance of the DPIM and 
 comparison with 7 harmonics HBFEM results. The FRCs are obtained with $\kappa_1=0.03\,\mum/\mus^2$ and  $\kappa_2=20,40\,\mum/\mus^2$, i.e., $\varepsilon_1=3.94\cdot 10^{-4}$ and  $\varepsilon_2=0.263,0.526$.
 (a): FRCs for the normalised modal coordinate of the master mode. The DPIM solution provides excellent results after introducing at least an order 3 approximation of the non-autonomous part. 
 (b-c): zoom near the peaks. Increasing the order improves the accuracy close to the FRCs peaks.
 (d-e): FRCs of the normalised modal coordinate of the slave mode. The DPIM approach allows reproducing the slave mode dynamics with a proper parametrisation of the non-autonomous part. The two forcing values are represented in separate subfigures for the sake of clearness.}
	\label{fig:multi_force}
\end{figure}

The external forcing amplitude is set equal to
$\kappa_1=1.5\,\mum/\mus^2$ i.e.  $\epsilon_1=0.019$, considering that $\max(\|\phivec_{B1}\|)=\phi=0.0898\mum$.
This load induces a midspan displacement of $1.08\mum$, which is about 1/6 of the arch thickness $H$. This vibration amplitude is sufficient to correctly excite the 2:1 superharmonic, which is also known to give rise to much smaller vibration amplitudes as compared to the primary resonance.
Fig.~\ref{fig:arch_12super} collects the FRCs considering different parametrisation orders, for a complex normal form style and reduction to a single NNM.
The FRCs plot the modal coordinate $u_{B1}=\phivec_{B1}^T \M \U$, normalised by the eigenmode amplitude $\phi$, and the characteristic length $H$.
The FRCs display how the proposed method converges towards the HBFEM results after adding the order 2 contribution in the non-autonomous forcing. 

To check the convergence trend of the method, different parametrisation orders have been considered. As expected,  $\varepsilon^1$ truncation fails to reproduce the super-harmonic effect, while $\varepsilon^2$ captures the main features. Increasing the order from $z^3$ to $z^5$  allows recovering the same nonlinearity content of the HBFEM solution. Further increasing the non-autonomous order to $\varepsilon^4$ provides a minor improvement that can be appreciated only at the FRC peak.
As in the previous example, the HBFEM convergence issues limit the effectiveness of the comparison between the two solutions close to the FRC peak. The ROMs solution is also in this case computed with MATCONT with the same set of parameters as in previous sections. For the sake of brevity, stability information has not been reported in this case.

An order $\mathcal{O}(z^5,\varepsilon^4)$ 
parametrisation needs less than 40 seconds to be computed and tracing the corresponding FRCs requires seconds. On the other hand, the 9 harmonics HBFEM solution takes approximately 3-5 minutes for each frequency value along the FRCs which results in 8-10 hours for a complete curve for this specific case. 

\subsection{Multimodal forcing of a shallow arch}
\label{sec:2F}

As a last example, the same arch structure as in the previous section is considered. It is now actuated with a more complex forcing shape involving two eigenmodes. The goal is to underline the effect of a strong forcing on a slave mode, when interested in the primary resonance of the master mode.
The body load has a spatial distribution proportional to a linear combination of bending mode B1 and the second symmetric bending mode B3. 
The corresponding FE model reads:
\begin{equation}
    \M \ddot{\U} + \C\dot{\U} +\K \U + \G(\U,\U) + \Hv(\U,\U,\U) = \kappa_1 \M \phivec_{B1} \cos(\omega t) + \kappa_2 \M \phivec_{B3} \cos(\omega t),
\end{equation}
This example corresponds to the one critically analysed in~\cite{opreni22high} with a $\varepsilon^1$ order truncation. 
The forcing factor $\kappa_2$ will be here further increased to enhance the observed discrepancies and show how the new method cures the limitations. Primary resonance of the fundamental mode is studied, such that $\omega \simeq \omega_{B1}$.

As in the previous examples, we introduce a normalisation of the external forcing to better highlight its magnitude with respect to the problem scale:
 \begin{equation}
    \epsilon_1=\frac{\phi_1}{H}\frac{\kappa_1}{\omega_{B1}^2},
\quad \mbox{with} \quad \phi_1=\max(\|\phivec_{B1}\|),
\end{equation}
 \begin{equation}
    \epsilon_2=\frac{\phi_3}{H}\frac{\kappa_2}{\omega_{B1}^2},
\quad \mbox{with} \quad \phi_3=\max(\|\phivec_{B3}\|).
\end{equation}
In the upcoming example, we will consider $\kappa_1=0.03$, that corresponds to $\epsilon_1=3.94 \cdot 10^{-4}$ and $\kappa_2=(20, 40)$  which gives $\epsilon_2=(0.263, 0.526)$ with $\phi_1=0.0898\mum$ and $\phi_3=0.0829\mum$. These loading scenarios lead to a peak midspan displacement of 0.981$\mum$ and 0.89$\mum$ respectively.
The load applied on the slave mode $\phivec_{B3}$
is strong enough to significantly alter the dynamics of the master mode. 
As a consequence, the maximum vibration amplitude decreases with increasing $\kappa_2$, such that the case $(\kappa_1,\kappa_2)=(0.03,20)$ corresponds to the upper curve with largest amplitude in Fig.~\ref{fig:multi_force}(a), while the lower curve corresponds to the loading scenario $(\kappa_1,\kappa_2)=(0.03,40)$. As in the previous examples, a mass-proportional Rayleigh damping model is assumed with $\alpha=\omega_{B1}/500, \beta=0$.

The results from the DPIM are compared with the ones achieved with an HB approach applied to the full finite element model (HBFEM) and 7 harmonics.
The FRCs are reported in Fig.~\ref{fig:multi_force}, where again different orders of truncation are considered for a complex normal form style parametrisation. To better highlight the difference between the proposed formulation and the $\varepsilon^1$ order truncation used for example  in~\cite{opreni22high}, we present 
the FRCs of the normalised modal coordinates of both the master mode $u_{B1}$ in Figs.~\ref{fig:multi_force}(a-c), and of the slave mode $u_{B3}$ in Figs.~\ref{fig:multi_force}(d-e).

The $\mathcal{O}(z^9,\varepsilon^0)$ solution is the same for the two forcing values $\kappa_2=(20,40) \mum/\mus^2$ for the projection on the master modal coordinate $u_{B1}$, thus resulting in a single curve shown in Fig.~\ref{fig:multi_force}(a); but is not for the slave modal coordinate $u_{B3}$, as shown in Figs.~\ref{fig:multi_force}(d,e). Indeed, as noted in Section~\ref{subsec:implementation}, the notation $\mathcal{O}(z^9,\varepsilon^0)$  refers to the case where the autonomous part is developed up to order~9, and the non-autonomous one is developed only up to order~1 in the non-autonomous variables $z^+,z^-$. This means that the component of the forcing along the slave mode only appears as a rigid translation of the manifold and it is therefore visible only on the slave mode, so the projection on the master mode is independent of the variation of the values of $\kappa_2$. From Figs.~\ref{fig:multi_force}(a,d,e), it appears clearly that this truncation is too crude and offers a prediction that is unacceptable as compared to the reference HBFEM solution.  
Note also that even the $\mathcal{O}(z^9,\varepsilon^1)$ solution, already used in~\cite{opreni22high}, begins to depart from the full-order solution at $\kappa_2=20 \mum/\mus^2$, and is then too far and unacceptable for $\kappa_2=40 \mum/\mus^2$. A minimal order to retrieve the correct result is given by $\mathcal{O}(z^9,\varepsilon^2)$, and convergence of the ROM is obtained for $\mathcal{O}(z^9,\varepsilon^3)$. This specific case thus underlines the gain brought by considering higher $\varepsilon$ order truncation in the non-autonomous part, even for cases corresponding to a primary resonance, with a strong excitation of the slave mode.

\section{Conclusion}

In this contribution, the parametrisation method for invariant manifold has been formulated in order to take into account the non-autonomous terms in a different manner as compared to previous developments shown {\it e.g.} in~\cite{BreunungHaller18,opreni22high}, where a first-order $\varepsilon$ development was imposed to deal with the forcing terms. Here, this assumption is bypassed using an automatic way of dealing with higher orders in $\varepsilon$ in the context of the parametrisation method. In particular, a dummy variable accounting for the exponential forcing term is included in the parametrisation algorithm as additional normal coordinates, and later replaced by its original definition to retrieve the time dependence of both invariant manifold parametrisation and reduced dynamics. In this way, the same routines used in the autonomous case can be readily extended to the non-autonomous one.

As a direct consequence of this new formulation, the method is not limited anymore to an $\varepsilon$ order development on the forcing amplitude, and is able to cope with superharmonic resonance. Indeed, the resonance relationship appearing from the solving of the homological equations makes appear multiple occurrences of the forcing frequency, a feature that was not present with a $\varepsilon^1$ order truncation.

The method has been implemented in a new version of the code \texttt{MORFE} which has been completely rewritten by also using the monomial representation with $\alphavec (p,k)$ vectors that were not used in~\cite{vizza21high,opreni22high}. In this new version of the code, the truncation can now be realised by selecting the order of the asymptotic expansion as well as the order of the $\varepsilon$ truncation.

Numerical results have been produced to demonstrate the gain brought about by this formulation as compared to previous developments shown in~\cite{vizza21high,opreni22high}. Consequently, the numerical developments have focused on illustrating that the method can deal with 3:1 and 2:1 superharmonic resonances respectively in a straight beam and a curved arch. Finally, an example with a two-mode shape forcing, already used in~\cite{opreni22high} and where the limits of the  $\varepsilon^1$ order truncation begin to appear, has been selected and pushed further in amplitude, definitely showing how the new implementation is now robust to a large range of forcing amplitudes.

The only limitation of the method seems now to be related to the fact that asymptotic expansions based on either graph or normal form styles have boundaries related to the assumed smallness of the normal variable. These boundaries are for example explored in~\cite{LamarqueUP,Stoychev:failing}. Another limitation for broad application to different nonlinearities is linked to the polynomial representation. Nevertheless, the method as it stands proposes very efficient ROMs for geometrically nonlinear structures, directly computable from the FE discretisation, able to converge to the desired accuracy thanks to arbitrary order expansion. Moreover, it directly computes the invariant manifold where the solutions are, thus offering the framework for exact reduced-order solutions.

\section*{Acknowledgments}
Giorgio Gobat is supported by the Joint Research Platform "Sensor sysTEms with Advanced Materials" (STEAM) - Politecnico di Milano and  STMicroelectronics S.r.l. Aur\'elien Grolet is thanked for numerous discussions and advice on the notations used.   Alessio Colombo is also thanked for his help in reading the preliminary drafts of the paper.

\section*{Funding}
The work received no additional funding.

\section*{Conflict of interest} 
The authors declare that they have no conflict of interest.

\section*{Data availability statement}
The implementation of the proposed method and the data required to run the numerical analyses presented in the paper are available from the repository~\cite{morfe}.
\section*{Author contributions}
Conceptualization: Alessandra Vizzacaro (lead), Giorgio Gobat (equal),  Attilio Frangi (equal), Cyril Touz\'e (lead); \\
Methodology: Alessandra Vizzacaro (lead), Giorgio Gobat (lead),  Attilio Frangi (lead), Cyril Touz\'e (lead); \\
Formal analysis and investigation: Alessandra Vizzacaro (equal), Giorgio Gobat (lead);  \\
Writing - original draft preparation: Alessandra Vizzacaro (equal), Giorgio Gobat (lead),  Attilio Frangi (equal), Cyril Touz\'e (lead); \\
Writing - review and editing: Attilio Frangi (lead), Cyril Touz\'e (lead);\\ 
Supervision:  Attilio Frangi (lead), Cyril Touz\'e (lead);

\bibliographystyle{unsrt}
\bibliography{main}

\appendix

\section{Comparisons on the methods to compute non-autonomous invariant manifolds}
\label{app:nonautreat}

The main purpose of this appendix is to compare more accurately the method proposed in the main text, to the more classical treatment of the non-autonomous term that relies on adding an auxiliary variable. For the sake of self-consistency, some of the equations already given in the main text are rewritten.  The development is again led on a first-order non-autonomous system to make the presentation simpler:
\begin{equation}
\label{eq:original_system}
    \B \dot{\y} = \A\y +\Q(\y,\y) 
    + \varepsilon \C \e^{\tilde{\lambda} t},
\end{equation}
with a single external forcing term, but its extension is straightforward in the case of multiple ones.

The most classical method to make autonomous a non-autonomous dynamical system consists of expanding the size of the phase space and considering an augmented state variable that will report on the forcing, as shown in all textbooks, see {\it e.g.}~\cite{gucken83,Wiggins} for instance. In the context of  a vibratory system, adding the forcing as an additional degree of freedom has been proposed in~\cite{ShawShock} to derive invariant-based ROMs. However, the procedure does not seem to have been further used by the authors, who then proposed another treatment in~\cite{JIANG2005H}, by adding a phase variable to make the system autonomous. 
In this appendix, we highlight the difference between this classical method and our proposal. To do so, the method proposed in~\cite{ShawShock} is first detailed. The first step is to make the original system autonomous and the second one is to compute the invariant manifold. 
Since a graph style parametrisation was considered in~\cite{ShawShock}, the realm is also extended here to the more general case of the parametrisation method that can also handle normal form style. At the end of the Appendix,  similarities and differences are finally discussed.

In order to make the system in~\eqref{eq:original_system} autonomous, one needs to introduce the following additional variable:
\begin{equation}
    \tilde{p} \doteq e^{\tilde{\lambda} t},
\end{equation}
and subsequently eliminate the explicit dependence on time by writing:
\begin{equation}
    \dot{\tilde{p}} = \tilde{\lambda}  \tilde{p},
\end{equation}
in conjunction with the initial condition: $\tilde{p}(0) = 1$. 
The new system reads:
\begin{equation}
    \begin{bmatrix}
        \B & \0\\
        \0 & 1
    \end{bmatrix}
    \begin{bmatrix}
        \dot{\y}\\\dot{\tilde{p}}
    \end{bmatrix}
    =
    \begin{bmatrix}
        \A & \varepsilon \C \\
        \0 & \tilde{\lambda}
    \end{bmatrix}
    \begin{bmatrix}
        {\y}\\\tilde{p}
    \end{bmatrix}
    +
    \begin{bmatrix}
        \Q(\y,\y) \\
        0 
    \end{bmatrix}.
    \label{eq:nonaut_made_aut}
\end{equation}

The parametrisation algorithm could now be applied directly to this autonomous system.
To seek a parametrisation of the invariant manifold from the space composed by both $\y$ and $\tilde{p}$, a nonlinear mapping should be introduced as
\begin{equation}
    \begin{bmatrix}
        {\y}\\\tilde{p}
    \end{bmatrix}
    =
    \begin{bmatrix}
        {{\W}({\z})}\\{\tilde{w}({\z})},
    \end{bmatrix},
    \label{eq:mappings_yp}
\end{equation}
where the last line concerned with the added coordinate $\tilde{p}$ has been isolated to keep track of it during the development. The dynamics on the sought manifold is expressed as usual with respect to the normal coordinates ${\z}$ as:
\begin{equation}
    \dot{{\z}} = {\f}({\z}).
\end{equation}

We will now demonstrate that if we apply the classical parametrisation algorithm in this autonomous setting, we will arrive at the same results that we reach with the strategy chosen in the paper. This is however only an algorithmic equivalence in terms of numerical computation, but not in terms of mathematical formulation, as discussed at the end of this Appendix. 

To apply the parametrisation algorithm to the system of Eq.~\eqref{eq:nonaut_made_aut}, the first step is to compute the linear master subspace. If $\tilde{\lambda}$ is resonant with an eigenvalue of the original system, then  a Jordan block in the eigenvalue matrix appears. The most interesting case and the one we will focus on here is the case of primary resonance. Say that an eigenvalue of the autonomous system, that we will denote with $\bar{\lambda}$, is equal or very close to the forcing eigenvalue $\tilde{\lambda}$, so that $\tilde{\lambda} \approx \bar{\lambda}$. Also, let us call $\bar{\Y}$ and  $\bar{\X}$ the right and left eigenvector of the original autonomous system associated to $\bar{\lambda}$.

One can see that the solution to the original eigenproblem associated with $\A$ and $\B$, is also a solution to the eigenproblem associated with Eq.~\eqref{eq:nonaut_made_aut}. In fact, the following equation is verified:
\begin{equation}
    \begin{bmatrix}
        \B & \0\\
        \0 & 1
    \end{bmatrix}
    \begin{bmatrix}
        \bar{\Y}\\0
    \end{bmatrix}
    \bar{\lambda}
    =
    \begin{bmatrix}
        \A & \varepsilon \C \\
        \0 & \tilde{\lambda}
    \end{bmatrix}
    \begin{bmatrix}
        \bar{\Y}\\0
    \end{bmatrix}.
    \label{eq:solution_trivial}
\end{equation}
Similarly $\begin{bmatrix}\bar{\X}^\star&
0\end{bmatrix}$ would be a solution of the left eigenproblem.

We now want to find the eigenvector associated with the eigenvalue $\tilde{\lambda}$. Let us call $\tilde{\Y}$ the upper part of this eigenvector related to the original variable $\y$, and $\tilde{y}$ the last entry of the eigenvector related to the added variable $\tilde{p}$. Also, let us remind that we are treating the case of $\tilde{\lambda} \approx \bar{\lambda}$. This is not a case of repeated eigenvalues \textit{per se} since the equality is only approximately fulfilled, but, in this context, near-resonances need to be treated in the same way exact resonances are, if we want to obtain meaningful results. This means that we are in the case of repeated eigenvalues with algebraic multiplicity equal to two. In such a case, one has to check the geometric multiplicity of the eigenvalues. If it is one, the eigenvalues admit a diagonal form. Otherwise, a Jordan block will appear in the eigenvalues matrix and the eigenvectors $[\bar{\Y}^T \; 0]^T$ and $[\tilde{\Y}^T\; \tilde{y}]^T$ will be coupled by an extra diagonal unitary term in the Jordan block. The geometric multiplicity of the eigenvalue can be simply checked and it is one if $\bar{\X}^\star\C = 0$, and two otherwise. If $\bar{\X}^\star\C = 0$, the method can be applied normally as one would do to square section beams, for instance. However, the condition $\bar{\X}^\star\C = 0$ means that the shape of the forcing is orthogonal to the mode it is supposed to force. Hence no direct spatial forcing is at hand, which is a very uncommon case. We are thus interested in the occurrence $\bar{\X}^\star\C \neq 0$, where the geometric multiplicity is two. A Jordan block then appears and the eigenproblem reads:
\begin{equation}
    \begin{bmatrix}
        \B & \0\\
        \0 & 1
    \end{bmatrix}
    \begin{bmatrix}
        \bar{\Y} & \tilde{\Y}\\
        0 & \tilde{y}
    \end{bmatrix}
    \begin{bmatrix}
        \bar{\lambda} & 1\\
        0 & \tilde{\lambda}
    \end{bmatrix}
    =
    \begin{bmatrix}
        \A & \varepsilon \C \\
        \0 & \tilde{\lambda}
    \end{bmatrix}
    \begin{bmatrix}
        \bar{\Y} & \tilde{\Y}\\
        0 & \tilde{y}
    \end{bmatrix}.
\end{equation}

The first column of the resulting system coincides with the system of Eq.~\eqref{eq:solution_trivial}, so it is automatically verified. The bottom row of the second column is also automatically verified by $\tilde{y}\tilde{\lambda} = \tilde{\lambda}\tilde{y}$. What we need to solve to obtain $\tilde{\Y}$ and $\tilde{y}$ is the top row of the second column, which gives the relationship:
\begin{equation}
    \B \bar{\Y} + \B \tilde{\Y} \tilde{\lambda}
    =
    \A \tilde{\Y} + \varepsilon \C \tilde{y}.
\end{equation}

However, this is an underdetermined system, so one last equation is needed. 
Since the orthogonality between eigenvectors should still be respected, we can complement the underdetermined system with an additional equation enforcing the orthogonality of the sought eigenvector $[\tilde{\Y}^T \; \tilde{y}]^T$ with respect to $\begin{bmatrix}\bar{\X}^\star&
0\end{bmatrix}$. Coupling the two equations together leads to:
\begin{equation}
    \begin{bmatrix}
        \B \tilde{\lambda}-\A & -\varepsilon \C
        \\
        \bar{\X}^\star\B & 0
    \end{bmatrix}
    \begin{bmatrix}
        \tilde{\Y} \\ \tilde{y}
    \end{bmatrix}
    =
    \begin{bmatrix}
        -\B \bar{\Y} 
        \\
        0
    \end{bmatrix}.
\end{equation}
It is possible to see that the upper part of the sought eigenvector $\tilde{\Y}$ is $\B$-orthogonal to $\bar{\X}$ and takes on the part of the forcing vector $\C$ which is also $\B$-orthogonal to $\bar{\X}$, whereas the lower part of the sought eigenvector is equal to:
\begin{equation}
    \tilde{y} = \dfrac{ \bar{\X}^\star \B \bar{\Y} }
    {\varepsilon \bar{\X}^\star \C }.
\end{equation}

This expression can be further simplified by noticing two things. First of all, the numerator is the arbitrary normalisation of $\bar{\Y}$ which we can set to 1. Secondly, from the definition of $\varepsilon$, the amplitude of $\C$ along the forced mode is also unitary, because the amplitude of the forcing must be quantified by $\varepsilon$. This means that the numerator and denominator simplify, leading to:
\begin{equation}
    \tilde{y} = \dfrac{1}
    {\varepsilon}.
\end{equation}

Now that we found the solution to the eigenproblem associated with Eq.~\eqref{eq:nonaut_made_aut}, we are finally able to see the effects that this choice of linear master subspace would have on the successive developments of the parametrisation algorithm. For the sake of readability, let us split the normal coordinates vector ${\z}$ into its two components as $\bar{z}$ the one related to the eigenvalue $\bar{\lambda}$ and $\tilde{z}$ the one related to the eigenvalue $\tilde{\lambda}$. Recalling now Eq.~\eqref{eq:mappings_yp}, we can see that the linear part of the mapping reads:
\begin{equation}
    \begin{bmatrix}
        {\y}\\\tilde{p}
    \end{bmatrix}
    =
    \begin{bmatrix}
        \bar{\Y}\\0
    \end{bmatrix}\bar{z}
    +
    \begin{bmatrix}
        \tilde{\Y}\\1/\varepsilon
    \end{bmatrix}\tilde{z}
    + \begin{bmatrix}
    \mathcal{O}(||{\z}||^2)
    \\
    \mathcal{O}(||{\z}||^2)
    \end{bmatrix},
\end{equation}

and the reduced dynamics reads:
\begin{equation}
    \begin{bmatrix}
        \dot{\bar{z}}\\\dot{\tilde{z}}
    \end{bmatrix}
    = 
    \begin{bmatrix}
        \bar{\lambda}\\0
    \end{bmatrix}
    \bar{z} 
    +
    \begin{bmatrix}
        1\\\tilde{\lambda}
    \end{bmatrix}
    \tilde{z} 
    + 
    \begin{bmatrix}
    \mathcal{O}(||{\z}||^2)
    \\
    \mathcal{O}(||{\z}||^2)
    \end{bmatrix}.
\end{equation}

Lastly, we must notice that setting to zero the higher order terms in the bottom rows of these last two equations solves the added equation for $\tilde{p}$, i.e. $\dot{\tilde{p}} = \tilde{\lambda}\tilde{p}$. This means that both the mapping $\tilde{w}({\z})$ and the reduced dynamics for the non-autonomous normal variable $\tilde{z}$ reduce to their linear part:
\begin{align}
& \tilde{p} = \dfrac{1}
    {\varepsilon}\tilde{z},
    \\
    & \dot{\tilde{z}}= 
        \tilde{\lambda}\tilde{z}.
\end{align}

Exactly like the strategy adopted in the paper, there is no  need to look for higher order terms in these two equations because they are automatically verified by their linear part, so only the mapping ${\W}({\z})$ and reduced dynamics for the normal variable $\bar{z}$ must be computed.
Moreover, since the original equation for $\tilde{p}$ is derived with an additional initial condition $\tilde{p}(0)=1$, this automatically implies that:
\begin{equation}
    \tilde{z}(0) = \varepsilon,
\end{equation}
which is also equivalent to what we assume in the strategy adopted in the paper, as $\tilde{z}=\varepsilon e^{\tilde{\lambda}t}$.

It is possible to demonstrate that pushing the developments to higher orders for the remaining functions ${\W}({\z})$ and the dynamics of $\bar{z}$, would also lead to the same numerical results as the strategy used in this paper. As a matter of fact, from an algorithmic point of view, the two techniques are fully equivalent. 

Having shown the numerical equivalence of the two strategies, we now want to discuss the key conceptual difference between them. Let us start by pointing out that, even though one could apply the parametrisation algorithm to any system, an actual convergence to a true invariant manifold is only proven under certain assumptions. There are indeed two main issues with the approach proposed in~\cite{ShawShock}. The first one is that, unless the forcing is decaying ($\text{Re}[\tilde{\lambda}]<0$), the theoretical results of the parametrisation method do not hold. In fact, these results lay on the assumption that all the eigenvalues of the system around the fixed point are stable. This is a rather uninteresting case as typically the forcing one $\tilde{\lambda}$ is purely imaginary. Secondly, since we are dealing with a local method, the results are only guaranteed in the vicinity of the fixed point, thus the initial condition with unitary amplitude ($\tilde{p}(0)=1$) is problematic. These issues are well known and the interested reader is referred to \cite{Haller2016,Llave2019} for more extensive discussions. 
As opposed to the strategy proposed in~\cite{ShawShock} where $\tilde{z}$ is a legitimate normal variable, we only introduce $\tilde{z}$ as a dummy variable in the algorithmic implementation and then retrieve the time-dependency of the parametrisation by substituting back its original definition. In this way, the condition on the stability of the fixed point required by the parametrisation method to ensure existence of an invariant manifold is not compromised. At the same time, no enforcement of a particular initial condition on the original system is imposed as we do not perform any modification to the original system.

\section{Detailed calculation for the homological equation of order $p$}
\label{app:pseudocode}

In this Appendix, we show how to deal with the generic order $p$ homological equation, at the level of the arbitrary monomial $(p,k)$. In particular, starting from Eq.~\eqref{eq:homologicpo1}, we want to isolate the single monomial by making all the terms in the equation explicit. Here, the complete derivation and all calculation details to go from Eq.~\eqref{eq:homologicpo1} to Eq.~\eqref{eq:zegoodhomolgicopo1:at} of the main text, are thus given for the sake of completeness.

The terms on the right-hand side of Eq.~\eqref{eq:homologicpo1} can be
rewritten as
\begin{subequations}\label{eq:zidevAQ}
    \begin{align}
        \A \left[\W (\z) \right]_p &= \sum_{k=1}^{m_p} \A \W^{(p,k)} \z^{\alphavec (p,k)}, \\
        \left[ \Q (\W,\W) \right]_p &= \sum_{k=1}^{m_p} \Q^{(p,k)} \z^{\alphavec (p,k)}.
    \end{align}
\end{subequations}
The quadratic terms are constructed from the product of lower-order terms, 
by exploiting the general relationship:
\begin{align}
\label{eq:codeQ}
    \Q \left(\W^{(p_1,k_1)} \z^{\alphavec (p_1,k_1)} \; , \;   \W^{(p_2,k_2)} \z^{\alphavec (p_2,k_2)}\right) = \Q \left(\W^{(p_1,k_1)} ,\W^{(p_2,k_2)}\right)  \z^{\alphavec (p_1,k_1) + \alphavec (p_2,k_2)},
\end{align}
for arbitrary orders $p_1$ and $p_2$. 
For a given order $p$, the $\Q^{(p,k)}$ terms can be computed in parallel for all $k$, according to:
\begin{align}
\label{eq:Q:at}
    \Q^{(p,k)} & = 
    \sum_{p_1=1}^{p-1} 
    \sum_{k_1,k_2=1}^{m_{p_1},m_{p_2}} 
    \Q \left(\W^{(p_1,k_1)} ,\W^{(p_2,k_2)}\right),  
    \\     
    & p_2: \quad p_2=p-p_1,
    \nonumber 
    \\
    & k : \quad \alphavec(p,k)=\alphavec (p_1,k_1) + \alphavec (p_2,k_2).
    \nonumber
\end{align}

Let us now consider the left-hand side term in Eq.~\eqref{eq:homologicpo1}, which is expanded to separate known and unknown terms as:
    \begin{align}
        \Bigl[ \nabla_{\z} & \W (\z) \f (\z) \Bigr]_p = \sum_{s=1}^{d+1} \Bigl[\frac{\partial \W}{\partial z_s} f_s (\z) \Bigr]_p
    \label{eq:deve000}\\
        &= \sum_{s=1}^{d+1} \Bigl[\left( \W^{(1,s)} + \frac{\partial \P[\dblnk]{\W(\z)}}{\partial z_s}  +  \frac{\partial \P[p]{\W(\z)}}{\partial z_s} \right)\left( \sum_{j=1}^{d+1} f_s^{(1,j)} z_j  +  \P[\dblnk]{f_s(\z)} +  \P[p]{f_s(\z)} \right)\Bigr]_p, \label{eq:nasty01}
        \nonumber
    \end{align}
where the shortcut notation  $\P[\dblnk]{\cdot}$ has been introduced to denote terms of order strictly larger than 1 and smaller than $p$. This three terms separation is needed to distinguish the known terms of order smaller than $p$ from the unknowns of order $p$. Since the operator $[\cdot]_p$ solely selects order $p$, only three terms from the product will stay instead of nine, the other ones being of lower orders, such that one can rewrite:
\begin{equation}\label{eq:ordreppp}
     \left[ \nabla_{\z} \W (\z) \f (\z) \right]_p = \sum_{s=1}^{d+1} \left( \W^{(1,s)} \P[p]{f_s(\z)} + \frac{\partial \P[p]{\W(\z)}}{\partial z_s} \left( \sum_{j=1}^{d+1} f_s^{(1,j)} z_j  \right)  + \left[ \frac{\partial \P[\dblnk]{\W(\z)}}{\partial z_s}    \P[\dblnk]{f_s(\z)}\right]_p \right).
\end{equation}
In this last equation, only the first two terms contain the unknowns. The last term is known as it involves products of lower orders. Let us separate the three terms of Eq.~\eqref{eq:ordreppp} as
\begin{equation}\label{eq:ordrepppN}
     \left[ \nabla_{\z} \W (\z) \f (\z) \right]_p = \N_1 (\z) + \N_2 (\z) + \N_3 (\z).
\end{equation}
One can rewrite explicitly
\begin{equation}
\label{eq:N3}
    \N_3 (\z)= \sum_{s=1}^{d+1} \sum_{p_W,p_f=1}^{o} \sum_{k_W,k_f=1}^{m_{p_W},m_{p_f}} \alpha_s (p_W,k_W) \W^{(p_W,k_W)} f_s^{(p_f,k_f)} \left[\z^{\alphavec (p_W,k_W)  + \alphavec (p_f,k_f) - \vectors{e}_s }  \right]_p,
\end{equation}
where the derivative 
\begin{equation}\label{eq:derivzsv2}
    \frac{\partial \z^{\alphavec (p_W,k_W)} }{\partial z_s} = \alpha_s (p_W,k_W) \z^{\alphavec (p_W,k_W) - \vectors{e}_{s}},
\end{equation}
has been used. The notation $\alphavec (p_W,k_W) - \vectors{e}_{s}$ is introduced to express the fact that the vector $\alphavec (p_W,k_W)$ has to be decreased by a unit value at the $s$-th component. In~\eqref{eq:derivzsv2}, $\alpha_s (p_W,k_W)$ is specified in order to keep track that this exponent comes from the derivative of the nonlinear mapping.  Since the operator $[\cdot]_p$ selects only order $p$, it is clear that from the double summations on all the possible orders $(p_W,p_f)$, only those fulfilling the relationship
\begin{equation}\label{eq:ppwpfv2}
    p=p_W+p_f-1,
\end{equation}
will be selected. Besides, the term $\N_3(\z)$ will produce $(d+1)$ monomials of order $p$, relative to the $s=1,\hdots,d+1$ index $k(s)$ verifying
\begin{equation}\label{eq:ksppwv2}
    \alphavec (p,k(s)) = \alphavec (p_W,k_W) + \alphavec (p_f,k_f) - \vectors{e}_s.
\end{equation}

\noindent
Consequently, the $\N_3$ term can be expressed as
\begin{equation}
\label{eq:R1v2}
\N_3 (\z)= \sum_{k=1}^{m_p}\N_3^{(p,k)} \z^{\alphavec(p,k)},
\end{equation}
where the coefficients $\N_3^{(p,k)}$ can be computed directly from Eq.~\eqref{eq:N3}
at a given order $p$ and for all $k$ at the same time:
\begin{align}
\label{eq:N3:at}
    \N_3^{(p,k)} & = 
    \sum_{s=1}^{d+1} 
    \sum_{p_W=2}^{p-1}
    \sum_{k_W,k_f=1}^{m_{p_W},m_{p_f}} 
    \alpha_s (p_W,k_W) \W^{(p_W,k_W)} f_s^{(p_f,k_f)},
    \\     
    & p_f: \quad p_f=p+1-p_W ,
    \nonumber 
    \\
    & k: \quad \alphavec(p,k)=\alphavec (p_W,k_W) + \alphavec (p_f,k_f) - \vectors{e}_s.
    \nonumber
\end{align}

The term $\N_1(\z)$ in Eqs.~\eqref{eq:ordreppp}-\eqref{eq:ordrepppN} can be rearranged by using the fact that the last line of the reduced dynamics, $f_{d+1}$, is linear by assumption and $[f_{d+1}(\z)]_p = 0$, for $p\geq2$. Also, as shown in Eq.~\eqref{eq:Wsolorder1Y}, the first-order linear mapping vectors $\W^{(1,s)}$ are the right eigenvectors $\Y_s$. Consequently, the term $\N_1$ can be rewritten as
\begin{equation}
    \N_1 (\z) = \sum_{k=1}^{m_p} \left(\sum_{s=1}^d \Y_s f_s^{(p,k)}\right) \z^{\alphavec (p,k)}.
\end{equation}
This term involves the unknown coefficients $f^{(p,k)}_s$ and will thus contribute to the left-hand side of the homological equation that gathers all unknown quantities.

The $\N_2 (\z)$ term can also be expanded, using the derivative given in Eq.~\eqref{eq:derivzsv2}, as:
\begin{equation}
\label{eq:N2dev:at}
    \N_2 (\z)= \sum_{s=1}^{d+1} \sum_{j=1}^{d+1} \sum_{k_W=1}^{m_p} \alpha_s (p,k_W) \W^{(p,k_W)} f_s^{(1,j)} \z^{\alphavec (p,k_W) - \vectors{e}_s + \vectors{e}_j}.
\end{equation}
The matrix $\f^{(1)}$ given in Eq.~\eqref{eq:f1dev00} is composed of a main diagonal collecting all eigenvalues $\lambda_1$ to $\lambda_d$ plus the forcing value $\tilde{\lambda}$ in $(d+1)\times (d+1)$ position. The only non-diagonal terms  appear in the last column with non-zero entries in $\f^{(1,d+1)}$ when a primary resonance occurs. This special structure allows rewriting Eq.~\eqref{eq:N2dev:at} as:
\begin{equation}
\label{eq:N2dev2:at}
    \N_2 (\z) = \sum_{k=1}^{m_p} \sum_{s=1}^{d+1} \alpha_s (p,k) \lambda_s \W^{(p,k)} \z^{\alphavec (p,k)} + \sum_{k_W=1}^{m_p} \sum_{s=1}^{d}  \alpha_s (p,k_W) \W^{(p,k_W)}  f_s^{(1,d+1)}  \z^{\alphavec (p,k_W)+\vectors{e}_{d+1} - \vectors{e}_s}. 
\end{equation}
In this last equation, the first term has been rewritten using the diagonal entries of $\f^{(1)}$ (i.e.\ $s=j$), while the non-zero terms in the last column $\f^{(1,d+1)}$ have been gathered in the second term. 

Let us denote as $\sigma^{(p,k)}$ the term:
\begin{equation}\label{eq:defsigmapk:app}
    \sigma^{(p,k)} = \sum_{s=1}^{d+1} \alpha_s \lambda_s .
\end{equation}

The first term in Eq.~\eqref{eq:N2dev2:at} is an unknown term that will hence contribute to the  left-hand side of the homological equation for a given monomial of order $p$ and index $k$.
The second term in Eq.~\eqref{eq:N2dev2:at} involves an order $p$ as well
so that $\N_2$ can be rewritten as:
\begin{equation}
\label{eq:R2v2:at}
    \N_2 (\z) = \sum_{k=1}^{m_p} \sigma^{(p,k)} \W^{(p,k)} \z^{\alphavec(p,k)}
    +
    \sum_{k=1}^{m_p} \N_2^{(p,k)} \z^{\alphavec(p,k)}.
\end{equation}
Unlike the coefficients $\N_3^{(p,k)}$ where all $k$ can be computed in parallel, see Eq.~\eqref{eq:N3:at},
 one needs to proceed here with a suitable ordering of the $\alphavec(p,k)$ vectors having the same order $p$
and, for a given $(p,k)$, the coefficients $\N_2^{(p,k)}$ are computed thanks to
\begin{subequations}\label{eq:N2pk:at}
\begin{align}
\N_2^{(p,k)} = 
\sum_{s=1}^{d}  \alpha_s (p,k_W) \W^{(p,k_W)}  f_s^{(1,d+1)}, \label{eq:N2pk:ata}
\\
 k_W: \quad \alphavec(p,k_W)=\alphavec (p,k)- \vectors{e}_{d+1} +\vectors{e}_s\,.
\label{eq:N2pk:at2b}    
\end{align}
\end{subequations}
Fortunately,
the homological equation for the $(p,k)$ combination
only depends on $\W^{(p,k_W)}$ such that 
Eq.~\eqref{eq:N2pk:at2b} holds.
By suitably ordering the list of $m_p$ vectors it can be guaranteed that 
the homological equation for $(p,k_W)$ is solved before 
addressing the case $(p,k)$. 
As a matter of fact, this result exploits the upper triangular nature of the linear reduced dynamics which makes the system here solvable iteratively. 
We can now write the order-$p$ homological equation at the level of an arbitrary monomial defined by $\alphavec (p,k)$. 
Let us denote as $\R^{(p,k)}$ the right-hand side term gathering all the known quantities as
\begin{equation}\label{eq:defRpkRHS}
    \R^{(p,k)} = \Q^{(p,k)} - \B \left(\N_3^{(p,k)} + \N_2^{(p,k)} \right).    
\end{equation}

Finally, the homological equation of order $p$ at the level of an arbitrary monomial reads
\begin{equation}
\label{eq:zegoodhomolgicopo1:app}
    \left( \sigma^{(p,k)} \B - \A \right) \W^{(p,k)} + \sum_{s=1}^d \B \Y_s f_s^{(p,k)} = \R^{(p,k)}.
\end{equation}

This equation is underdetermined and considerations on its solvability are reported in the main text in Section~\ref{sec:orderp}.

\section{Considerations on the order of truncation of the normal variables and the dummy variables}
\label{app:truncation}

In this Appendix we want to give more insight on the problem of the treatment of the different asymptotics linked to the amplitude of the dummy variables and of the normal variables, and how this is embedded in the strategies used to solve out the problem. From a mathematical point of view, the smallness of the forcing is automatically embedded in the fact that the asymptotic expansion is performed around $\z=\0$, with $\z$ now including both $\bar{\z}$ and $\tilde{\z}$. In order to treat $\bar{\z}$ (related to the dynamics) and $\tilde{\z}$ (related to the forcing) with the same assumptions in terms of expansions, the only thing left to define is the relationship between the smallness of $\varepsilon=||\tilde{\z}||$ as compared to that of ${\rho}=||\bar{\z}||$. Typically in nonlinear vibrations, large amplitude displacement can be observed even with very small forcing level if the forcing frequency resonate with one of the system's natural frequencies. For this reason, it makes sense to define $\varepsilon$ to be smaller than ${\rho}$ in the asymptotic sense. Let us write such relationship in the general case as:
\begin{equation}
    \mathcal{O}(\varepsilon) = \mathcal{O}({\rho}^{m}),
    \label{eq:orders}
\end{equation}
with $m$ a nonzero integer that can be tuned. 

Notice that the choice of $m$ affects the truncation order for the normal variables. In fact, it follows from Eq.~\eqref{eq:orders}, that a truncation at order ${o}$ for ${\rho}$ corresponds to:
\begin{equation}
    \mathcal{O}({\rho}^{{o}})
    =\mathcal{O}(\varepsilon^{{o}/m})
    =
    \mathcal{O}({\rho}^{\bar{p}}\varepsilon^{(o-\bar{p})/m}), \qquad \forall \bar{p}=0,\dots,o
\end{equation}

Basically, given $\bar{p}$ being the order of the autonomous variables, and $\tilde{p}$ that of the forcing variables, only the monomials respecting $\bar{p}+m\tilde{p}\leq o$, are kept in the expansion. An example of which monomials will be left in the expansion under this approach is given in Tab.~\ref{tab:truncationa}-\ref{tab:truncationb}. From an asymptotic point of view, this is the most rigorous approach to the choice of the truncation order of the monomials. However, it is not intuitive to foresee what value of $m$ is appropriate for a given case. In fact, fixing $\bar{o}$ and $m$ automatically decides $\tilde{o}$, the maximum order of the $\varepsilon$ expansion. This is in opposition to what one would like to do, which is to freely choose $\bar{o}$ and $\tilde{o}$. In fact, $\tilde{o}$ has a much clearer physical meaning than $m$ in most cases, for instance, if one is interested in the third superharmonic resonance of a system, at least a third order in the forcing expansion is needed. It is then easier to set the maximum needed order for the forcing than selecting a value for  $m$. For this reason, the approach we use in the code is to set $m=1$ in all cases and to add another constraint to the maximum order of the epsilon expansion: $\tilde{p}\leq\tilde{o}$. An example of the monomials that would be left by following this approach is given in Tab.~\ref{tab:truncationc}. For the sake of completeness, it is worth mentioning that a third approach could be also adopted, which is related to the idea that there is no relation between the smallness of $\varepsilon$ and that of $\rho$. In such a case, only constraints on the single variables are given ($\bar{p}\leq\bar{o}$ and $\tilde{p}\leq\tilde{o}$) without enforcing any constraint on their relationship. This is what we call the disjoint approach and an example of it is given in Tab.~\ref{tab:truncationd}.


\begin{table}[ht]
\centering
\subfloat[Asymptotic approach with $m=3$: $\bar{p}+3\tilde{p}\leq o$]{\label{tab:truncationa}
    \begin{tabular}{|c|c|c|c|}
    \hline
    $\mathcal{O}(\rho^0\varepsilon^0)$
    & 
    $\mathcal{O}(\rho^0\varepsilon^1)$ 
    & 
    $\mathcal{O}(\rho^0\varepsilon^2)$
    & 
    \phantom{$\mathcal{O}(\rho^0\varepsilon^3)$}
    \\
    \hline
    $\mathcal{O}(\rho^1\varepsilon^0)$& $\mathcal{O}(\rho^1\varepsilon^1)$& &
    \\
    \hline
    $\mathcal{O}(\rho^2\varepsilon^0)$& $\mathcal{O}(\rho^2\varepsilon^1)$& &
    \\
    \hline
    $\mathcal{O}(\rho^3\varepsilon^0)$& $\mathcal{O}(\rho^3\varepsilon^1)$& &
    \\
    \hline
    $\mathcal{O}(\rho^4\varepsilon^0)$& & &
    \\
    \hline
    $\mathcal{O}(\rho^5\varepsilon^0)$& & & 
    \\
    \hline
    $\mathcal{O}(\rho^6\varepsilon^0)$& & &
    \\\hline
    \end{tabular}
    }
\hfill
\subfloat[Asymptotic approach with $m=2$: $\bar{p}+2\tilde{p}\leq o$]{\label{tab:truncationb}
    \begin{tabular}{|c|c|c|c|}
    \hline
    $\mathcal{O}(\rho^0\varepsilon^0)$
    & 
    $\mathcal{O}(\rho^0\varepsilon^1)$ 
    & 
    $\mathcal{O}(\rho^0\varepsilon^2)$
    & 
    $\mathcal{O}(\rho^0\varepsilon^3)$
    \\
    \hline
    $\mathcal{O}(\rho^1\varepsilon^0)$& $\mathcal{O}(\rho^1\varepsilon^1)$& $\mathcal{O}(\rho^1\varepsilon^2)$&
    \\
    \hline
    $\mathcal{O}(\rho^2\varepsilon^0)$& $\mathcal{O}(\rho^2\varepsilon^1)$& $\mathcal{O}(\rho^2\varepsilon^2)$&
    \\
    \hline
    $\mathcal{O}(\rho^3\varepsilon^0)$& $\mathcal{O}(\rho^3\varepsilon^1)$& &
    \\
    \hline
    $\mathcal{O}(\rho^4\varepsilon^0)$& $\mathcal{O}(\rho^4\varepsilon^1)$& &
    \\
    \hline
    $\mathcal{O}(\rho^5\varepsilon^0)$& & & 
    \\
    \hline
    $\mathcal{O}(\rho^6\varepsilon^0)$& & &
    \\\hline
    \end{tabular}
    }
\hfill
\subfloat[The approach we use: $\bar{p}+\tilde{p}\leq o \;\&\; \tilde{p} \leq  \tilde{o}$]{\label{tab:truncationc}
    \begin{tabular}{|c|c|c|c|}
    \hline
    $\mathcal{O}(\rho^0\varepsilon^0)$
    & 
    $\mathcal{O}(\rho^0\varepsilon^1)$ 
    & 
    $\mathcal{O}(\rho^0\varepsilon^2)$
    & 
    $\mathcal{O}(\rho^0\varepsilon^3)$
    \\
    \hline
    $\mathcal{O}(\rho^1\varepsilon^0)$& $\mathcal{O}(\rho^1\varepsilon^1)$& $\mathcal{O}(\rho^1\varepsilon^2)$& $\mathcal{O}(\rho^1\varepsilon^3)$
    \\
    \hline
    $\mathcal{O}(\rho^2\varepsilon^0)$& $\mathcal{O}(\rho^2\varepsilon^1)$& $\mathcal{O}(\rho^2\varepsilon^2)$& $\mathcal{O}(\rho^2\varepsilon^3)$
    \\
    \hline
    $\mathcal{O}(\rho^3\varepsilon^0)$& $\mathcal{O}(\rho^3\varepsilon^1)$& $\mathcal{O}(\rho^3\varepsilon^2)$& $\mathcal{O}(\rho^3\varepsilon^3)$
    \\
    \hline
    $\mathcal{O}(\rho^4\varepsilon^0)$& $\mathcal{O}(\rho^4\varepsilon^1)$& $\mathcal{O}(\rho^4\varepsilon^2)$& 
    \\
    \hline
    $\mathcal{O}(\rho^5\varepsilon^0)$& $\mathcal{O}(\rho^5\varepsilon^1)$& & 
    \\
    \hline
    $\mathcal{O}(\rho^6\varepsilon^0)$& & & 
    \\\hline
    \end{tabular}
    }
\hfill
\subfloat[Disjoint approach: $\bar{p}\leq o \;\&\; \tilde{p} \leq  \tilde{o}$]{\label{tab:truncationd}
    \begin{tabular}{|c|c|c|c|}
    \hline
    $\mathcal{O}(\rho^0\varepsilon^0)$
    & 
    $\mathcal{O}(\rho^0\varepsilon^1)$ 
    & 
    $\mathcal{O}(\rho^0\varepsilon^2)$
    & 
    $\mathcal{O}(\rho^0\varepsilon^3)$
    \\
    \hline
    $\mathcal{O}(\rho^1\varepsilon^0)$& $\mathcal{O}(\rho^1\varepsilon^1)$& $\mathcal{O}(\rho^1\varepsilon^2)$& $\mathcal{O}(\rho^1\varepsilon^3)$
    \\
    \hline
    $\mathcal{O}(\rho^2\varepsilon^0)$& $\mathcal{O}(\rho^2\varepsilon^1)$& $\mathcal{O}(\rho^2\varepsilon^2)$& $\mathcal{O}(\rho^2\varepsilon^3)$
    \\
    \hline
    $\mathcal{O}(\rho^3\varepsilon^0)$& $\mathcal{O}(\rho^3\varepsilon^1)$& $\mathcal{O}(\rho^3\varepsilon^2)$& $\mathcal{O}(\rho^3\varepsilon^3)$
    \\
    \hline
    $\mathcal{O}(\rho^4\varepsilon^0)$& $\mathcal{O}(\rho^4\varepsilon^1)$& $\mathcal{O}(\rho^4\varepsilon^2)$& $\mathcal{O}(\rho^4\varepsilon^3)$
    \\
    \hline
    $\mathcal{O}(\rho^5\varepsilon^0)$& $\mathcal{O}(\rho^5\varepsilon^1)$& $\mathcal{O}(\rho^5\varepsilon^2)$& $\mathcal{O}(\rho^5\varepsilon^3)$
    \\
    \hline
    $\mathcal{O}(\rho^6\varepsilon^0)$& $\mathcal{O}(\rho^6\varepsilon^1)$& $\mathcal{O}(\rho^6\varepsilon^2)$& $\mathcal{O}(\rho^6\varepsilon^3)$
    \\\hline
    \end{tabular}
    }
\caption{Visualisation of the order of the monomials kept in the expansion for different truncation approaches. Case of $o=6$, $\bar{o}=6$, $\tilde{o}=3$.}\label{tab:truncation}
\end{table}

\section{Treatment of the homological equations in the modal space}
\label{app:projmodalhomo}

In this appendix, more details are given to the solution of the homological equations by using a projection onto the modal space. Only the case of first-order systems is detailed.

To begin with, the solution to Eq.~\eqref{eq:order1supcoldp1}  written in the modal space for a better understanding of the appearance of resonance relationship, is first highlighted. 
Let us denote as $\xivec^{(1,d+1)}$ the projection of the unknown mapping vector  $\W^{(1,d+1)}$ in the modal space:
\begin{equation}
        \xivec^{(1,d+1)} = \X_{tot}^\star \B \W^{(1,d+1)}, 
\end{equation}
where $\X_{tot}$ collects the $D$ left eigenvectors, and has been defined in Eq.~\eqref{eq:defmatrixofeigenvectX}.
Let us also introduce  the modal forcing vectors $\F$ as
\begin{equation}
        \F =  \X_{tot}^\star \C, \\
\end{equation}
Eq.~\eqref{eq:order1supcoldp1} can be projected onto the modal space by left-multiplication with the matrix  $\X_{tot}^\star$, yielding:
\begin{equation}\label{eq:order1supcolproj}
    \X_{tot}^\star\left( \tilde{\lambda} \B - \A   \right) \W^{(1,d+1)} = \X_{tot}^\star\C - \X_{tot}^\star\B \Y\f^{(1,d+1)}.
\end{equation}
Using the relationship defining the left eigenvectors, Eq.~\eqref{eq:lefteigendef}, one can write
\begin{equation}
    \X_{tot}^\star\A \W^{(1,d+1)} = \Lambdavec \xivec^{(1,d+1)},
\end{equation}
where $\Lambdavec$ is the diagonal matrix containing the $D$ eigenvalues $\{ \lambda_1, \hdots, \lambda_D \}$ as defined in Eq.~\eqref{eq:defBiglambda}. Assuming that the modes have been normalised with respect to $\B$ such that the matrix $\D$ appearing in Eq.~\eqref{eq:normalizedFOS} is the identity matrix, Eq.~\eqref{eq:order1supcolproj} finally reads:
\begin{equation}\label{eq:order1modal}
    \left[ \tilde{\lambda} \I_D - \Lambdavec \right] \xivec^{(1,d+1)} = \begin{bmatrix}
        F_1 - f_1^{(1,d+1)} \\ \vdots \\ F_d - f_d^{(1,d+1)} \\ F_{d+1} \\ \vdots \\ F_{D}
    \end{bmatrix}.
\end{equation}
This equation is interpreted in the same lines as general solutions at order $p$ for the homological equation, the only difference being that, with the new treatment of the forcing proposed, a first resonance already appears at first order. One can separate the contributions of the first $d$ lines of Eq.~\eqref{eq:order1modal} (the tangent space, relative to the master modes) to the last $D-d$ lines corresponding to the normal space relative to the slave modes. In the normal space, one can write $\forall \, j \in [d+1,D] $:
\begin{equation}
    (\tilde{\lambda} - \lambda_j ) \xi_j^{(1,d+1)} = F_j.
\end{equation}
By definition there is no primary resonance condition between the forcing frequency and the slave modes, consequently $\tilde{\lambda} \neq \lambda_j$, and thus the solution for the unknown mapping in the modal space is known. In the tangent space, one arrives at, $\forall \, j \in [1,d] $:
\begin{equation}\label{eq:zisolmodalo1}
    (\tilde{\lambda} - \lambda_j ) \xi_j^{(1,d+1)} = F_j - f_j^{(1,d+1)}.
\end{equation}
Two unknowns are still present, the solution is thus given by following the discussion on resonance: if there is no primary resonance, meaning that the forcing value $\tilde{\lambda}$ is far from any eigenvalue of the master modes, then Eq.~\eqref{eq:zisolmodalo1} is solved by setting $f_j^{(1,d+1)}$ to zero and solving for $\xi_j^{(1,d+1)}$. In this case, the linear term of the reduced dynamics becomes diagonal, meaning that the forcing is not an order 1 term for the dynamics onto the manifold. The other case is the resonant one. Let us denote as $\mR^{(1,d+1)}$ the set of of primary resonant modes, {\it i.e.} 
\begin{equation}
    \mR^{(1,d+1)} = \left\{ r_1, \hdots, r_l \right\}, \quad \mbox{such that} \quad \lambda_{r_j} = \tilde{\lambda} ,
\end{equation}
with $\Y_{r_j}$ the associated resonant right eigenvectors. The cardinality of $\mR^{(1,d+1)}$ is here $l$ which can be larger than 1 in order to take into account the particular case of degenerate modes where the eigenvalues has multiplicity $l$. In the notation $\mR^{(1,d+1)}$, the upperscript $(1,d+1)$ refers to the fact that order-1 term is at hand (first entry), and the resonance relationship is here considered with $\tilde{\lambda}$, which is the eigenvalue of the $d+1$ linear monomial, $z_{d+1} = \tilde{z}$. This notation naturally generalises to any monomial of arbitrary order.

The solution to Eq.~\eqref{eq:zisolmodalo1} is then given by:
\begin{align}
    \mbox{if}\; r \in \mR^{(1,d+1)}, & \quad \xi_r^{(1,d+1)}=0, \quad \mbox{and} \quad f_r^{(1,d+1)} = F_r,\\
    \mbox{if}\; r \notin \mR^{(1,d+1)}, & \quad \xi_r^{(1,d+1)}=\frac{F_r}{\tilde{\lambda} - \lambda_r}, \quad \mbox{and} \quad f_r^{(1,d+1)} = 0.
\end{align}
Note that for the resonant modes for $r \in \mR^{(1,d+1)}$, the condition $\xi_r^{(1,d+1)}=0$ translates back to the physical space to
\begin{equation}
    \X_r^\star \B \W^{(1,d+1)} = 0,
\end{equation}
which is the condition used to augment the singular system and make it solvable while imposing this vanishing condition in the modal space.

The same projection can be used to shed light on the solution of the homological equation of order $p$, written at the level of an arbitrary monomial defined by $\alphavec (p,k)$, Eq.~\eqref{eq:zegoodhomolgicopo1:at}. Left-multiplying Eq.~\eqref{eq:zegoodhomolgicopo1:at} by $\X_{tot}^\star$ for projection, it reads:
\begin{equation}\label{eq:zegoodhomolgicopproj}
    \sigma^{(p,k)}  \X_{tot}^\star\B \W^{(p,k)}  - \X_{tot}^\star\A \W^{(p,k)} + \sum_{s=1}^d  \X_{tot}^\star\B \Y_s f_s^{(p,k)} =  \X_{tot}^\star\R^{(p,k)}.
\end{equation}
Defining $\Svec^{(p,k)} = \X_{tot}^\star\R^{(p,k)}$ as the modal projection of the left-hand side, and using the orthonormality relationships, one arrives at
\begin{equation}\label{eq:zegoodhomolgicopproj2}
    \left(\sigma^{(p,k)} \I_D - \Lambdavec   \right) \xivec^{(p,k)}  + \begin{bmatrix}
        f_1^{(p,k)} \\ \vdots \\ f_d^{(p,k)} \\ 0\\ \vdots \\ 0
    \end{bmatrix} = \Svec^{(p,k)}.
\end{equation}
The solution of this equation follows then the classical discussion related to the use of the parametrisation method for invariant manifolds~\cite{Haro}. In the normal space, {\it i.e.} $\forall s=d+1,\hdots,D$, the equation writes simply as
\begin{equation}\label{eq:orderpnormalsp}
    (\sigma^{(p,k)} - \lambda_s) \xi_s^{(p,k)} = S_s^{(p,k)}.
\end{equation}
It is assumed that no cross-resonance exists, such that the relationship $\lambda_s \simeq \sigma^{(p,k)}$ is never fulfilled for the slave modes. Then Eq.~\eqref{eq:orderpnormalsp} has a single solution. The case is different in the tangent space relative to the master coordinates since we will have more unknowns than equations. Indeed, $\forall s=1,\hdots,d$, one has
\begin{equation}\label{eq:orderptangentsp}
    (\sigma^{(p,k)} - \lambda_s) \xi_s^{(p,k)}  + f_s^{(p,k)}= S_s^{(p,k)}.
\end{equation}
Let us introduce $\mR^{(p,k)}$ the resonant set, which collects all the $r$ indexes such that the nonlinear resonance relationship $\lambda_r \simeq \sigma^{(p,k)}$ is fulfilled:
\begin{equation}\label{eq:defRpk}
    \mR^{(p,k)} = \{ r \in [1,d] \; | \;  \lambda_r \simeq \sigma^{(p,k)} \}.
\end{equation}
The two main styles of parametrisations follow the guidelines given below. In a graph style parametrisation, then Eq.~\eqref{eq:orderptangentsp} is solved by setting
\begin{equation}
    \forall s=1,\hdots,d,\quad \xi_s^{(p,k)} =0, \quad f_s^{(p,k)}= S_s^{(p,k)}.
\end{equation}
This means that, in the graph style parametrisation, the set $\mR_{graph}^{(p,k)}$  is filled with all indexes of the master modes, irrespective of the fact that the resonance relationship is effectively fulfilled
\begin{equation}
    \mR_{graph}^{(p,k)} = \{ 1, \; 2, \; \hdots, \; d \}.
\end{equation}

In a normal form style parametrisation, Eq.~\eqref{eq:orderptangentsp} is solved following the dichotomy:
\begin{align}
    \forall \, s \, \in \mR^{(p,k)}, & \quad \xi_s^{(p,k)} =0, \quad f_s^{(p,k)}= S_s^{(p,k)};\\
    \forall \, s \, \notin \mR^{(p,k)}, & \quad f_s^{(p,k)}= 0, \quad \xi_s^{(p,k)} = \frac{S_s^{(p,k)}}{\sigma^{(p,k)} - \lambda_s}.
\end{align}
In the normal form style, $\mR^{(p,k)}$ strictly follows the definition \eqref{eq:defRpk} and gathers only the resonant monomials, with the view of offering the most simple reduced dynamics: the normal form.

As a final note, one can observe that imposing $\xi_s^{(p,k)} =0$ when $s \, \in \mR^{(p,k)}$, translates in the original coordinates as:
\begin{equation}
    \X_s^\star \B \W^{(p,k)} = 0,
\end{equation}
which is the orthogonality condition used to solve the augmented system in a direct approach.

\section{Second order mechanical problem: further simplifications using real modes}
\label{app:2ndorderreal}

The aim of this appendix is to show how one can further simplify the order-$p$ homological equations to be solved for mechanical systems, Eq.~\eqref{eq:mecho2orderpfinal}, where the property relating displacement and velocity has already been used to halve the size of the problems to be solved. By doing so, the link with the previous developments reported in~\cite{vizza21high,opreni22high} will be completely made explicit. Let us assume here that the mass, damping and stiffness matrices of the mechanical problem are symmetric such that
\begin{equation}
    \M^t = \M, \quad \C^t = \C, \quad \K^t = \K. 
\end{equation}
The real normal mode $\phivec_k$ with associated radian eigenfrequency $\omega_k$  are the solutions of the conservative eigenproblem
\begin{equation}
    (\K - \omega_k^2 \M)\phivec_k = \0.
\end{equation}
For the sake of simplicity, it is here assumed that the modes are mass-normalised, and the damping matrix is diagonalised by the real normal modes, such that:
\begin{equation}
    \phivec_k^t \M \phivec_l = \delta_{kl}, \quad \phivec_k^t \K \phivec_l = \omega_k^2\delta_{kl}, \quad \phivec_k^t \C \phivec_l = 2\xi_k \omega_k\delta_{kl},
\end{equation}
where the modal damping ratio $\xi_k$ has been introduced. The complex-conjugate eigenvalues of the damped problem come by pairs and read
\begin{equation}
    \lambda_k = -\xi_k \omega_k \pm i \omega_k \sqrt{1-\xi_k^2}.
\end{equation}
They are solutions of the problem
\begin{equation}
    \left( \lambda_k^2 \M + \lambda_k \C + \K \right) \phivec_k = \0.
\end{equation}
A last and important property shared by the modes has been demonstrated in Appendix~B of~\cite{vizza21high} and will be used to simplify the equations:
\begin{equation}\label{eq:zesimplifkitu}
    \left( \lambda_k \M + \C  \right) \phivec_k = -\bar{\lambda}_k \M \phivec_k.
\end{equation}

With these assumptions, the left and right eigenvectors $\X_k$ and $\Y_k$ used in the main text can be completely rewritten as function of the real normal modes of the mechanical problem. The matrix of right master eigenvectors $\Y$ introduced in Eq.~\eqref{eq:YrightMASTER} can be advantageously rewritten for mechanical systems. It can be first divided in two parts, by separating the first $N$ lines $\Y^V$ corresponding to the velocity, to the last $N$ lines $\Y^U$ corresponding to the displacement. Since a vibration mode corresponds to a pair of complex eigenvalues, selecting $n$ master modes $(\phivec_1,\hdots,\phivec_n)$ leads to a matrix $\Y$ having $2n$ columns, such that 
\begin{equation}
    \Y^U = \begin{bmatrix}
        \phivec_1 & \hdots & \phivec_n & \phivec_1 & \hdots & \phivec_n
    \end{bmatrix}.
\end{equation}
The left eigenvectors $\X_k$ can also be constructed from the mass-normalised real normal modes as it has been done in~\cite{vizza21high,opreni22high}. In this case, for the velocity part (first $N$ lines), one has
\begin{equation}
    \X_k^V = \frac{1}{\lambda_k - \bar{\lambda}_k} \phivec_k.
\end{equation}
Note that the scalar in front of the eigenmode shape just results from the fact that the left eigenvectors are in this case constructed from the mass-normalised mode. If one uses a normalisation with respect to the $\B$ matrix of the first-order problem, then this term will disappear. Since this is only a problem related to normalisation, it will not affect the next developments which can be conducted in the same manner. Since the goal here is to make the link with~\cite{vizza21high}, the real modes mass-normalised will be kept for this calculation.

Let us focus on the added term appearing in Eq.~\eqref{eq:mecho2orderpfinal}, ${\X^V_{\mR}}^\star \M \Y^U_{\mR} $, in order to show how this term simplifies to retrieve the problem formulation derived in~\cite{vizza21high,opreni22high}. The aim of the derivation here is to show that one can retrieve the formula to solve the second-order mechanical problem given in~\cite{vizza21high}. Since the augmented problem to be solved depends on the style of parametrisation selected and on the filling of the set $\mR$, one has to distinguish here the case of complex normal form, real normal form and graph style parametrisation.

Let us begin with the complex normal form style parametrisation, assuming a single master mode labeled $r$ for simplicity and no internal resonance, as in Section 4.3~of~\cite{vizza21high}. If the monomial under study $\alphavec(p,k)$ is resonant with $r$, then $\mR_{\text{CNF}} = \{ r \}$, and the term simplifies to
\begin{equation}
    {\X^V_{\mR}}^\star \M \Y^U_{\mR} = {\X^V_r}^t \M \Y^U_r = \frac{1}{\lambda_k - \bar{\lambda}_k} \phivec_r^t \M \phivec_r = \frac{1}{\lambda_k - \bar{\lambda}_k}.
\end{equation}
Besides, the upper term $\left[\left( \sigma+ \lambda_{r\in\mR}\right)\M + \C \right] \Y^U_{\mR}$ also simplifies using Eq.~\eqref{eq:zesimplifkitu}:
\begin{equation}
    \left[\left( \sigma+ \lambda_{r\in\mR}\right)\M + \C \right] \Y^U_{\mR} = \left[\left( \sigma+ \lambda_r \right)\M + \C \right] \Y^U_r = (\sigma - \bar{\lambda}_r) \M \phivec_r.
\end{equation}
The second line of Eq.~\eqref{eq:mecho2orderpfinal} makes appear the product by the left eigenvector ${\X^V_{\mR}}^\star $, so the factor $\frac{1}{\lambda_k - \bar{\lambda}_k}$ can be easily simplified, such that one is able to write the augmented system for the CNF as
\begin{equation}
\begin{bmatrix}
\sigma^2\M + \sigma \C + \K & \quad(\sigma-\bar{\lambda}_r)\M \phivec_r 
\\
(\sigma - \bar{\lambda}_r) \phivec_r^t\M & 1 
\end{bmatrix}
\begin{bmatrix}
        \Psivec^{(p,k)} \\
        f^{(p,k)}_{r} 
\end{bmatrix}
=
\begin{bmatrix}
         \Xivec^{(p,k)}\\
        \phivec_{r}^t  \M \muvec^{(p,k)}\\
\end{bmatrix},
\label{eq:homo_reduced_cnf}
\end{equation}
which is exactly the formulation given in~\cite{vizza21high}, see Eq.~(80). 

In the case of the real normal form style, then if the monomial considered $\alphavec (p,k)$ is resonant with the master mode $r$, then the choice $\mR_{\text{RNF}} = \{ r \,\, r^{\star} \}$ is made. Once again the terms appearing in Eq.~\eqref{eq:mecho2orderpfinal} can be easily simplified and one retrieves the formulation given in Eq.~(83) in~\cite{vizza21high}. For the graph style parametrisation, the filling of $\mR$ contains all the master modes and the supplementary terms of the augmented system makes finally appear a $2n\times 2n$ matrix composed of four $n \times n $ blocks with the identity matrix, as in Eq.~(88) of~\cite{vizza21high}.

\section{Normalisation based on structural theory}
\label{sec:an_beam}

The aim of this appendix is to justify the nondimensional parameter introduced in Eq.~\eqref{eq:epsload}, which is used in order to monitor the importance of the forcing term. Whereas the characteristic length used to monitor the importance of the transverse displacement with respect to geometric nonlinearities is generally given as the thickness (except for unconstrained beams in the axial direction where it stands as the length, see {\it e.g.}~\cite{ReviewROMGEOMNL} for a discussion), quantifying the importance of forcing with a nondimensional number is less standard. The goal here is thus to introduce the scaling by referring to the simple case of a beam of length $L$, restrained in the axial direction. In such case, the von K{\`a}rm{\`a}n-Mettler assumptions can be used to derive the equations of motion as
\begin{equation}
    \rho A \frac{\partial^2 v}{\partial t^2} + EJ \frac{\partial^4  v}{\partial x^4} - \frac{EA}{2L} \left[ \int_0^L \left( \frac{\partial v}{\partial x} \right)^2 \mathrm{d}x \right] \frac{\partial^2  v}{\partial x^2} = F(x,t).
\end{equation}
In this equation, $\rho$ stands for the density, $A$ the area of the cross-section, $E$ the Young modulus, $J$ the quadratic moment,  $v(x,t)$ is the transverse displacement and $F(x,t)$ the distributed external force. Note that the boundary conditions are not specified since they are not needed for the development shown here. Assuming a single mode expansion as $v(x,t)=\phi(x)q(t)$, with $\phi$ the non-dimensional shape function, one easily arrives to an oscillator equation
\begin{equation}\label{eq:duffbeam}
    M\ddot{q} + K q + K_3 q^3 = f,
\end{equation}
where
\begin{subequations}
    \begin{align}
        M &= \int_0^L \rho A \phi^2\,\text{d}x,\\
        K &= \int_0^L EJ \phi '''' \phi \text{d}x,\\
        K_3 &= -\frac{EA}{2L} \int_0^L \left( \phi' \right)^2  \text{d}x \int_0^L \phi \phi''  \text{d}x.
    \end{align}
\end{subequations}
Dividing by the mass, adding dissipation and assuming a harmonic forcing with excitation frequency $\omega$, Eq.~\eqref{eq:duffbeam} can be written as:
\begin{equation}
\ddot{q}+ 2\xi\omega_0 \dot{q} +\omega_0^2 q+ k_3 q^3 =\kappa \cos(\omega t),
\end{equation}
with $k_3=K_3/M$ and $\omega_0^2 = K/M$. 
Finally, let us introduce the nondimensional time $\tau=\omega_0 t$ to have a unitary linear stiffness, and let us make the displacement nondimensional with  $Q=q \max(\phi)/L_{CH}$, where $\max(\phi)$ is the maximal value of the eigenmode shape, and $L_{CH}$ the characteristic length (in general the thickness, 
but for unconstrained beams the length). The non-dimensional equation then reads:
\begin{equation}
Q''+ 2\xi Q' +Q+ \left(\frac{L_{CH}}{\max(\phi)}\right)^2\frac{k_3}{\omega_0^2} Q^3 =\frac{\max(\phi)}{L_{CH}}\frac{\kappa}{\omega_0^2}\cos(\tau).
\end{equation}
One can see that the nondimensional parameter appearing in front of the forcing amplitude is the one used in Eq.\eqref{eq:epsload}.

\end{document}